\documentclass[english]{article}
\usepackage[T1]{fontenc}
\usepackage[latin9]{inputenc}
\usepackage[letterpaper]{geometry}
\usepackage{float}
\usepackage{amsmath}
\usepackage{graphicx}
\usepackage{color}
\usepackage{epstopdf}
\usepackage{adjustbox}
\usepackage{diagbox}
\usepackage{multirow}
\usepackage{array}
\usepackage{mathrsfs}
\usepackage{amsmath, amssymb, amsthm}
\usepackage{subfigure}
\usepackage{tikz}

\graphicspath{{./Figures/}}

\usepackage{amsfonts} 
\usepackage{color}
\definecolor{black}{rgb}{0,0,0}

\definecolor{red}{rgb}{1,0,0}

\definecolor{blue}{rgb}{0,0,1}

\usepackage{multirow}

\makeatother

\usepackage{babel}
\usepackage{authblk}
\title{Online mixed multiscale finite element method with oversampling and its applications}

\author[1]{Yanfang Yang}
\author[2]{Shubin Fu \thanks{Corresponding Author}}
\author[2]{Eric T. Chung}
\affil[1]{School of Mathematics and Information Science, Guangzhou University, Guangzhou, People's Republic of China}
\affil[2]{Department of Mathematics, The Chinese University of Hong Kong, Hong Kong SAR}

\begin{document}
		\maketitle
\begin{abstract}
	In this paper, we consider an  online basis enrichment mixed generalized multiscale method with oversampling, for solving flow problems in highly heterogeneous porous media. This is an extension of the online mixed generalized multiscale method \cite{online_mixed}. The multiscale online basis functions
	are computed by solving a Neumann problem in an over-sampled domain, instead of a standard neighborhood of a coarse face. We are motivated by the restricted domain decomposition method.
	 Extensive numerical experiments are presented to demonstrate the
	performance of our methods for both steady-state flow, and two-phase flow and transport problems.
	
\end{abstract}
{{\it keywords: } multiscale; mixed finite element; oversampling; heterogeneous media.}	

\section{Introduction}
In many scientific and engineering applications, multiple scales and high contrast are common features. For instance, in gas and oil production, the reservoir properties, such as the permeability, can be detailed at multiple scales, varying from inches to miles. Many model reduction techniques, such as upscaling and multiscale methods, are developed for
the solution of  such kind of problems. These model reduction approaches aim to reduce the degrees of freedom and solve the problem on a coarse grid. For example, in upscaling methods \cite{durfolsky1991homo,wu2002analysis}, one averages the media properties based by some rules and then solve problems on a coarse grid. In multiscale methods \cite{egw10,efendiev2009multiscale,Arbogast_two_scale_04,chung2015mixed,chen2003mixed,jennylt03,Wheeler_mortar_MS_12, Arbogast_PWY_07, mortar_elliptic}, one still solves the problems on a coarse grid, with precomputed multiscale basis functions that are constructed locally on the fine grid and  carry the multiscale information of the media.
	
	In this paper, we  introduce an oversampled enrichment algorithm  in the framework of generalized multiscale finite element method (GMsFEM) in solving  flow problems in the mixed
	formulation in heterogeneous media. We first compute a coarse grid solution with offline basis functions  which are formed by following the GMsFEM. Then we iteratively compute basis functions using residuals from  the previous solution in the online stage, thus we call it an online method. The method in this paper is  an extension of the online mixed generalized multiscale method \cite{online_mixed}. Here we employ the oversampling technique for the computation of the online basis functions, and study the application
	of this method for practical two phase simulation problems.

	We follow a mixed finite element framework of the flow problem. An important motivation of using mixed methods is to ensure mass conservation, which is essential for flow problems. The main idea of mixed multiscale  finite element method is to construct multiscale basis functions for each coarse face supported in the neighborhood of this face.  Mixed GMsFEM is a generalization of the  classical mixed multiscale  finite element method \cite{aarnes04, chen2003mixed}. In the classical method, only one basis function per coarse face is used to capture the multiscale features of the media.  However, it is shown \cite{ge09_2, egh12} that for some cases such as non-separable scales and long channels, one basis function per coarse face is not sufficient to capture all the features. For this sake,  GMsFEM  \cite{chung2015mixed} is proposed which allows more basis functions per coarse face.
	
	The computation of basis functions in mixed GMsFEM consists of offline and online stages. In the offline stage, offline basis functions are formed by solving a series of local spectral problems. These functions carry the multiscale features and they  can be reused for other input parameters to solve the equation. In the online stage, online basis functions are constructed based on the parameters. However, as pointed out in \cite{online_cg,online_mixed,online_dg}, there are multiscale problems in which an offline process fails to yield an accurate representations of solutions. The reason for this is that offline computations are carried out locally. Therefore, various basis enrichment approaches are developed.
	
	In \cite{cel_13}, an adaptive algorithm is proposed to enrich the solution space by adding basis functions which are precomputed in the offline stage. In \cite{online_cg}, adaptive methods in the continuous Galerkin framework, using the residuals of the previous solutions to form new online basis functions are discussed. Results show that these methods accelerate the convergence rate of GMsFEM significantly.  Based on the idea of using residuals, there are also related approaches developed for the discontinuous Galerkin formulation in \cite{online_dg} and mixed formulation in \cite{online_mixed},  where	an online adaptive method to enrich the solution space is proposed. This online adaptive method forms new basis functions by projecting the previous solution on the space of divergence free functions. In this paper, we follow the idea of paper \cite{online_mixed}, except that we use an oversampling technique in the computation of the online basis functions to further accelerate the convergence rate.
	
	Oversampling is introduced for multiscale finite element methods in \cite{hw97}, and it is proven that oversampling can improve the accuracy of multiscale methods. Oversampling is introduced  in the context of GMsFEM in \cite{eglp13}, where oversampling is used in the offline stage to compute the snapshot and offline spaces. Other methods such as mixed GMsFEM and mixed mortar GMsFEM involving oversampling for offline computation can be found in \cite{chung2015mixed, mortar_elliptic}. In this paper, we  discuss the use of oversampling for online basis computation. This is motivated by the restricted domain decomposition method presented in \cite{cai1999restricted}, where the local preconditioners are obtained by utilizing larger domain than standard domain to perform computation and then taking restriction onto  standard domain. In this way, one can achieve faster convergence.
	
	For our online basis enrichment method with oversampling, a set of non-overlapping subdomains is selected. New basis functions are constructed on each of these subdomains. In particular, we solve a residual problem on an oversampled domain covering an coarse face. Then we restrict the normal trace of the solution to the coarse face. The new basis function is obtained by solving a Neumann problem with the normal trace as boundary condition for the coarse face on this coarse faces's neighborhood.

	We show numerical experiments for various heterogeneous permeability fields. We compare
	 online enrichment without oversampling and online enrichment with oversampling and show that online with oversampling has much faster convergence.  We study both high and low conductivity inclusions and channels in the domain. It is shown that we can get the same convergence rate for different order of high contrast orders. We also consider a 3D two-phase flow and transport problem with permeability field from the last 50 layers of the SPE10 benchmark model \cite{Aarnes2005257}. We demonstrate that by adding a few online basis functions for the initial problem, without updating solution space as time advancing, we can get accurate results.
	
	The rest of the paper is organized as follows. The basic idea of mixed GMsFEM  is presented in the next section, including the definition of coarse and fine grids, the construction of the snapshot space and the offline space. In section 3 we introduce the online basis enrichment method with oversampling. Extensive numerical results are given in section 4. We conclude the paper in section 5.

\section{Preliminaries}
We consider the following system of flow problem in a mixed formulation:
\begin{equation}\label{eq:orgional_equation}
\begin{split}
\kappa^{-1}v+\nabla p&=0 \quad \text{in } D,\\
\text{div}(v)&=f\quad \text{in } D,
\end{split}
\end{equation}
with the homogeneous Neumann boundary condition $v\cdot n=0$ on $\partial D$, where $\kappa$
 is a high-contrast permeability field, $D\subset\mathbb{R}^{d}(d=2,3)$  is the computational domain and $n$ is the unit outward normal
 vector of the boundary of $D$.

 The basic idea of Generalized multiscale finite element method(GMsFEM) is to construct multiscale basis functions on a local coarse region. To present our method, we first introduce fine and coarse grids.
Let $\mathcal{T}_H$ be a partition of $D$ into coarse blocks 
$K_i$ with diameter $H_i$ so that $\overline{D}=\cup_{i=1}^N\overline{K}_i$,
where $N$ is the number of coarse blocks. $\mathcal{T}_H$ is called the coarse grid, on which our coarse grid discretization will be defined.
We call $E_H$ a coarse
face of the coarse element $K_i$ if $E_H = \partial K_i \cap \partial K_j $ or $E_H= \partial K_i \cap \partial{D}$.
Let $\mathcal{E}_H(K_i)$ be the set of all coarse faces on the boundary of the coarse block $K_i$ and $\mathcal{E}_H=\cup_{i=1}^N\mathcal{E}_H(K_i)$
be the set of all coarse faces. For our mixed GMsFEM, the basis functions are constructed on $\omega_i$, which are the two coarse elements that share a common face, i.e.,
$$\omega_i=\cup \{K\in \mathcal{T}_H: E_i\in \partial K\}, i=1,2,\cdots, N_e,$$ where $N_e$ is the number of coarse faces. We also define an over-sampled domain of a coarse face $E_i$ as $\omega_i^{+}$. Note that $\omega_i^+$ is not necessarily a union of coarse blocks.

\begin{figure}[H]
	\centering
	\resizebox{0.7\textwidth}{!}{
		\begin{tikzpicture}[scale=0.7]
		\filldraw[fill=orange, draw=black] (2,6) rectangle (10,10);
		\filldraw[fill=blue, draw=black] (4,4) rectangle (8,12);
		
		\draw[step=4,black, line width=0.8mm] (0,0) grid (12,16);
		\draw[step=1,gray, thin] (0,0) grid (12,16);
		\draw[ultra thick, red](4, 8) -- (8,8);
		\draw [->,dashed, thick](12.4,10) -- (6.5, 8);
		\node at (16.8,10.2)  {\huge $E_{i}$: coarse edge (red)};

		\draw [->,dashed, thick](12.4,6) -- (6.5, 7.2);
		\node at (19.2,6.2)  {\huge $\omega_{i}$: Coarse neighborhood (blue)};
		
		\draw [->,dashed, thick](12.4,2.6) -- (8.5, 6.4);
		\node at (19.8,2.8)  {\huge $\omega_{i}^+$: Oversampling region (orange)};
		
		\node at (2.5,9.5)  {\huge$\omega_{i,1}^+$} ;
		\node at (2.5,6.5)  {\huge$\omega_{i,2}^+$} ;			
		\node at (6.5,10.5)  { \huge $K_i^1$} ;
		\node at (6.5,5.5)  { \huge $K_i^2$} ;
		\end{tikzpicture}
	}
	\caption{ Illustration of a coarse edge $E_i$, and its coarse neighborhood $\omega_{i}$, oversampling neighborhood $\omega_{i}^+$.}
	\label{grids}
\end{figure}
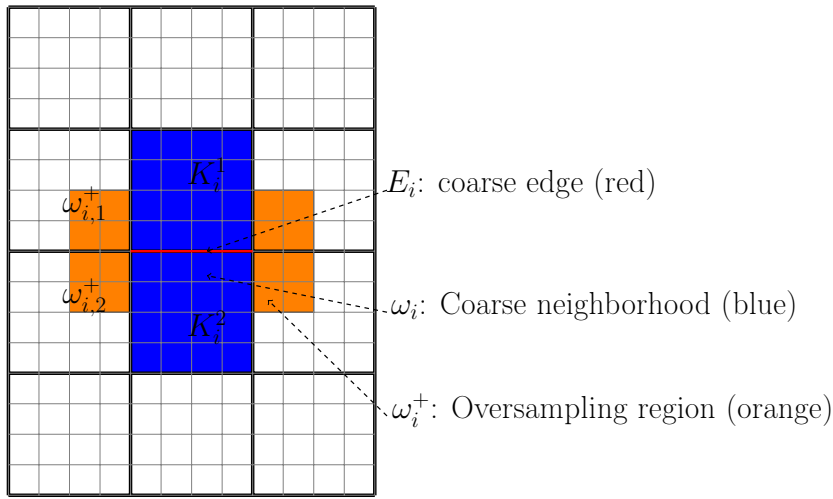

We further partition each coarse block $K_i$ using a finer mesh $\mathcal{T}_h(K_i)$ with mesh size $h_i$.
Let $\mathcal{T}_h=\cup_{i=1}^N\mathcal{T}_h(K_i)$ be the union of all these partitions, which is a fine mesh partition of the domain $D$.
We use $h=\text{max}_{1\leq i\leq n}h_i$ to denote the mesh size of $\mathcal{T}_h$.
In addition, we let
$\mathcal{E}_h(K_i)$
be the set of all faces of the partition $\mathcal{T}_h(K_i)$ and $\mathcal{E}_h^0(K_i)$ be the set of all interior faces of the partition
$\mathcal{T}_h(K_i)$ and let $\mathcal{E}_h=\cup_{i=1}^N\mathcal{E}_h(K_i)$ be the set of all faces in the partition $\mathcal{T}_h$. 
Figure~\ref{grids} gives an illustration of the constructions of the two grids for the case of 2D. The black lines represent the coarse grid, and the gray lines represent the fine grid.
For each coarse face $E_i$, we define a coarse neighborhood $\omega_i$ as the union of all coarse blocks having the face $E_i$.
Figure~\ref{grids} shows a coarse neighborhood $\omega_i$ in the blue  color.
The oversampled neighborhood $\omega_i^+$ (yellow region in Figure~\ref{grids}) is the joint of two subdomains $\omega_{i,1}^+$
and $\omega_{i,2}^+$ such that $E_i^+=\partial\omega_{i,1}^+\cap\partial\omega_{i,2}^+$.
$E_i^+$ is the union of fine faces that covers a coarse face $E_i$.

We let $V_h\times Q_h$ be the standard lowest-order Raviart-Thomas space
for the approximation of (\ref{eq:orgional_equation}) on the fine grid $\mathcal{T}_h$.
Then, the fine-grid solution $(v_h, p_h)$ satisfies
\begin{equation}
\begin{split}
\int_D\kappa^{-1}v_h\cdot w-\int_D\text{div}(w)p_h&=0, \quad\quad\quad\quad \forall w\in V_h^0,\\
\int_D\text{div}(v_h)q&=\int_Dfq, \quad\quad \forall q \in Q_h.
\end{split}
\end{equation}
where $v_h\cdot n=0$ on $\partial D$ and $V_h^0=V_h\cap\{v\in V_h: v\cdot n=0  \text{ on } \partial D\}$.
The above system can be written in terms of matrix representations as
\begin{equation}
\left[\begin{array}{cc}
A_h & B_h^{T} \\
B_h & 0\\

\end{array}\right]\left[\begin{array}{cc}v_{hr}\\p_{hr}\end{array}\right] =\left[\begin{array}{cc}0\\F_h\end{array}\right]
\label{eq:fine_system}
\end{equation}
where $v_{hr}$ and $p_{hr}$ are vectors of coefficients in the expansions
of the solutions $v_h$, $p_h$ in the space $V_h$ and $Q_h$.
We will use the fine-grid solution $(v_h, p_h)$ as the reference solution, and we will compare the accuracy
of the multiscale solution $(v_H, p_H)\in (V_{H}, Q_H)$ against the fine-grid solution. Here $V_{H}, Q_H$  are multiscale solution spaces for the pressure
$p$ and the velocity $v$  in Equation (\ref{eq:orgional_equation}). Next we give the definition of $V_{H}, Q_H$.

First $Q_H$  consists of functions which are constant on each coarse grid block.  The construction of the space $V_H$ follows the framework of the mixed GMsFEM. Generally, it has two steps.
First, we construct a snapshot space  $V_{\text{snap}}$ consisting of an extensive set of functions which are able  to approximate the solution $v$. To construct the snapshot space,
we construct a set of snapshot functions $\{\psi_{\text{snap},i}^j\}_{j=1}^{N_i}$ supported in the coarse grid neighborhood $\omega_i.$  We define the local snapshot
space for $\omega_i$ as $V^{i}_{\text{snap}}=\text{span}(\{\psi_{\text{snap},i}^j\}_{j=1}^{N_i}).$ We call the span of all the local snapshot functions as the snapshot space
 $V_{\text{snap}}=\oplus_{E_i \in \mathcal{E}_H }V^{i}_{\text{snap}}.$  The second step is to construct an offline space $V_H$ by reducing the large snapshot space to a smaller one.

With the pressure space $Q_H$ and the velocity space $V_H$, the mixed GMsFEM is to find  $v_{H}\in V_{H}$ and $p_{H}\in Q_H$ such that
\begin{equation}
\begin{split}
\int_D\kappa^{-1}v_{H}\cdot w-\int_D\text{div}(w)p_{H}&=0, \quad\quad\quad\quad \forall w\in V_{H}^0,\\
\int_D\text{div}(v_{H})q&=\int_Dfq, \quad\quad \forall q \in Q_H,
\end{split}
\end{equation}
with where $v_H\cdot n=0$ on $\partial D$ and $V_H^0=V_H\cap\{v\in V_H: v\cdot n=0  \text{ on } \partial D\}$.
The matrix form of the above coarse system is
\begin{equation}
\left[\begin{array}{cc}
R^TA_hR & R^TB_h^{T}G_H \\
G_H^TB_hR & 0\\

\end{array}\right]\left[\begin{array}{cc}v_{Hr}\\p_{Hr}\end{array}\right] =\left[\begin{array}{cc}0\\G_H^TF_h\end{array}\right]
\label{eq:coase_system}
\end{equation}
where $R$ stores all the multiscale basis functions,
$G_H$ is the restriction operator from $Q_H$ into $Q_h$.
$v_{Hr}$ and $p_{Hr}$ are vectors of coefficients in the expansions
of the solutions $v_H$, $p_H$ in the space $V_H$ and $Q_H$.

In the following sections, we present the construction of the snapshot space  $V_{\text{snap}}$ and the offline space  $V_{H}$ in detail.

\subsection{Snapshot space}\label{snapshot_space}
In this section, we discuss the construction of the snapshot space  $V_{\text{snap}}$ which consists of basis functions up to the resolution of the fine grid faces on the coarse grid faces.
We construct the local snapshot spaces $V^{i}_{\text{snap}}$ by solving a set of local problems on each coarse neighborhood  $\omega_i$, and then define the global snapshot space $V_{\text{snap}}=\oplus_{E_i \in \mathcal{E}_H }V^{i}_{\text{snap}}.$

Let $E_i \in {\cal E}_H$, which can be written as a union of fine-grid faces, i.e., $E_i=\cup_{l=1}^{J_i}e_l$, where $J_i$ is the total number of fine-grid faces on $E_i$ and $e_i$ represents a fine-grid face.
We solve the following problems
\begin{equation}
\begin{split}
\kappa^{-1} v_l^{i}+\nabla p_l^i&=0,\quad\quad \text{in} \quad \omega_i,\\
\text{div}(v_l^{i})&=\alpha_l^i,\quad \quad \text{in}\quad \omega_i.
\end{split}
\label{snapshot_eqn}
\end{equation}
subject to the homogeneous Neumann boundary condition $v_l^{i} \cdot n_i=0$ on $\partial \omega_i.$ The above problem is solved separately on each coarse-grid block contained in $\omega_i$, so that the snapshot basis consists of solutions
of local problems with all possible boundary conditions on the face $E_i$ up to the fine-grid resolution.  To solve the equation (\ref{snapshot_eqn}) on $K\subset \omega_i$, an additional boundary condition
$v_l^{i} \cdot n_i=\delta_l^i$ on $E_i$ is used, where $\delta_l^i$ is defined by
\begin{equation}
\delta_l^i=\left\{
\begin{aligned}
 1,\quad \text{on} \quad e_l,\\
0, \quad \text{on} \quad E_i \backslash e_l,\\
\end{aligned}
\quad l=1,2, \cdots, J_i,
\right.
\end{equation}
and $n_i$ is a fixed unit-normal vector for $E_i.$ The constant $\alpha_l^i$ in equation (\ref{snapshot_eqn}) is chosen to satisfy the compatible condition $\int_k \alpha_l^i=\int_{\partial K}v_l^i\cdot n_i.$
The solutions of the above local problems form the local snapshot space $V^{i}_{\text{snap}}$, from which we get  $V_{\text{snap}}=\oplus_{E_i \in \mathcal{E}_H }V^{i}_{\text{snap}}.$ Next, we discuss the
derivation of the offline space from $V_{\text{snap}}$.

\subsection{Offline space}\label{offline_space}
As we mentioned earlier, the snapshot space $V_{\text{snap}}$ is a large space with dimension comparable to the fine resolution. We will perform a dimension reduction on $V_{\text{snap}}$ to get a smaller space. This reduced space is called
the offline space.
The reduction is accomplished by solving a local spectral problem on each coarse grid neighborhood $\omega_i.$
The local spectral problems is to find real number $\lambda$ and $v\in V^{i}_{\text{snap}}$ such that
\begin{equation}
a(v,w)=\lambda s(v,w),\quad \forall w \in V^{i}_{\text{snap}},
\label{spectral_1}
\end{equation}
where $a(\cdot,\cdot)$ and $s(\cdot,\cdot)$ are symmetric positive definite bilinear operators on $V^{i}_{\text{snap}}\times V^{i}_{\text{snap}}.$
Specifically, we take
\begin{equation}
\begin{split}
a(v,w)&=\int_{E_i} \kappa^{-1}(v\cdot n_i)(w\cdot n_i),\\
s(v,w)&= \frac 1 {H}\left(\int_{\omega_i} \kappa^{-1}v\cdot w+\int _{\omega_i} \text{div}(v)\text{div}(w)\right),
\end{split}
\label{spectral_eqn}
\end{equation}
for $v, w\in V^{i}_{\text{snap}}$, and $n_i$ is the fixed unit normal vector for $E_i$.

After solving the spectral problem (\ref{spectral_1}) in $\omega_i$, we arrange the eigenvalues in  ascending order
$$\lambda^i_1< \lambda^i_2\cdots <\lambda^i_{J_i}.$$
Let $\phi^i_1, \phi^i_2,\cdots, \phi^i_{J_i}$ be the corresponding eigenfunctions. We select the first $l_i$ eigenfunctions to form the offline space $V^i_{\text{off}}$,
i.e., $V^i_{\text{off}}=\text{span}\{\phi^i_1, \phi^i_2,\cdots, \phi^i_{l_i}\}$.
The global offline space is  $V_{\text{off}}=\oplus_{E_i \in \mathcal{E}_H }V^{i}_{\text{off}}.$

Suppose the matrix form of the local snapshot space $V^{i}_{\text{snap}}$ is
$R_{\text{snap}}^i=[\Psi_{\text{snap}}^{i,1},\cdots,\Psi_{\text{snap}}^{i,J_i}]$.  If we follow the same procedure as in Section \ref{snapshot_space} on the domain $\omega_i^+,$ we can get the snapshot
space ${V}^{i,+}_{\text{snap}}$ with basis functions supported on $\omega_i^+$.
More specifically,
we solve the following problem on the oversampling neighborhood $\omega_i^+$
\begin{equation}
\begin{split}
\kappa^{-1} v_l^{i,+}+\nabla p_l^{i,+}&=0,\quad\quad \text{in} \quad \omega_i^+,\\
\text{div}(v_l^{i,+})&=\beta_l^i,\quad \quad \text{in}\quad \omega_i^+.
\end{split}
\label{snapshot_eqn1}
\end{equation}
subject to the boundary condition
\begin{equation}
v_l^{i,+} \cdot n_i=0,\quad \text{on} \quad\partial\omega_i^+.
\end{equation}
We solve the  problem above on $\omega_{i,1}^+$ and $\omega_{i,2}^+$ separately with additional boundary condition
$v_l^{i,+} \cdot n_i=\delta_l^i$ on $E_i^+$, where $E_i^+$ is the intersection of $\omega_{i,1}^+$ and $\omega_{i,2}^+$, see Figure~\ref{grids}. 
The constant $\beta_l^i$ is chosen to guarantee the compatible condition such that the equation (\ref{snapshot_eqn1}) is solvable.
Then the local oversampled snapshot space ${V}^{i,+}_{\text{snap}}$ can be formed with the solution
of the above local problem.
Denote $\widehat{V}^{i,+}_{\text{snap}}$ be the divergence free subspace of ${V}^{i,+}_{\text{snap}}$.
We denote  the matrix form of the oversampled local snapshot space ${V}^{i,+}_{\text{snap}}$ as
$R_{\text{snap}}^{i,+}=[\Psi_{\text{snap},1}^{i,+},\cdots,\Psi_{\text{snap},J_i^+}^{i,+}]$,
where $J_i^+$ is the number of fine-grid faces on $E_i^+$.

In the next section, we will give the construction of  a residual driven online space
$V_{H}$ with the oversampling techniques.

\section{Oversampling online iterative algorithm}
In this section, we discuss an enrichment algorithm for the velocity space which constructs new basis functions based on the solution of the previous enrichment level. These new basis functions are formed in the online computation stage, therefore we call them online basis. With the enrichment of the online basis functions, we can achieve a fast convergence. The novelty of this paper is that we construct the online basis by using an oversampling approach.

For each region $\Omega \subseteq D$, we denote  $V_{\Omega}$ as the space of functions in $V_{\text{snap}}$ which are supported in $\Omega$, i.e.,    $V_{\Omega}=\oplus_{\omega_i \subseteq \Omega }V^{i}_{\text{snap}}.$  Let $\widehat{V}_{\Omega}$ be the divergence free subspace of $V_{\Omega}$.
Then we define the linear functional $R_{\Omega}$ on $V_{\Omega}$ by
\begin{equation*}
R_{\Omega}(v)=\int_{\Omega}\kappa^{-1}v_{H}\cdot v-\int_{\Omega}\text{div}(v)p_{H}, \quad \forall v\in V_{\Omega}.
\end{equation*}
We note that if we restrict $R_{\Omega}$ on $\widehat{V}_{\Omega}$, then
\begin{equation*}
R_{\Omega}(v)=\int_{\Omega}\kappa^{-1}v_{H}\cdot v, \quad \forall v\in \widehat{V}_{\Omega}.
\end{equation*}

\textbf{Online iterative algorithm with oversampling:}
We start with iteration number $l=0$. $V_{H}^0$ is taken to be the offline space $V_{\text{off}}$ defined in Section \ref{offline_space}.

{\bf Step 1:} Find the multiscale solution in the current spaces $V_{H}^l$ and $Q_H$. That is to find  $v_{H}^l\in V_{H}^l$ and $p_{H}^l\in Q_H$ satisfying
 \begin{equation}
 \begin{split}
 \int_D\kappa^{-1}v_{H}^l\cdot w-\int_D\text{div}(w)p_{H}^l&=0, \quad\quad\quad\quad \forall w\in V_{H}^l,\\
 \int_D\text{div}(v_{H}^l)q&=\int_D fq, \quad\quad \forall q \in Q_H.
 \end{split}
 \end{equation}

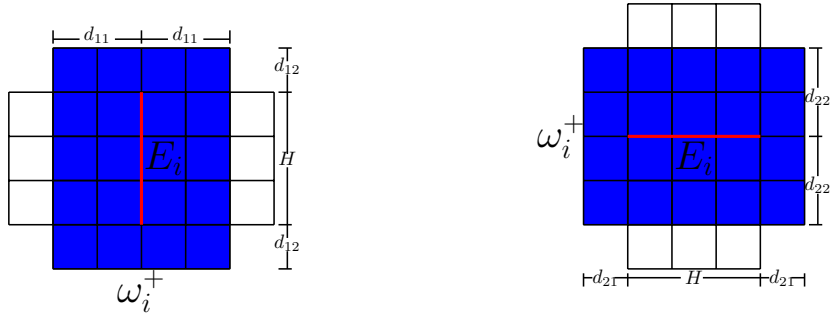
\begin{figure}[H]
		\centering
		\resizebox{0.7\textwidth}{!}{
			\begin{tikzpicture}[scale=0.8]
			\filldraw[fill=blue, draw=black] (1,0) rectangle (5,5);
			\draw[step=1,black, thick] (1,0) grid (5,5);
			\draw[step=1,black, thick] (0,1) grid (6,4);
			\draw[ultra thick, red](3, 1) -- (3,4);
			\node at (3.5,2.5)  { \huge $E_{i}$} ;
			\node at (12.5,3)  { \huge $\omega_{i}^+$} ;
			\draw[ thick, black](1, 5.1) -- (1,5.4);
			\draw[ thick, black](1, 5.2) -- (1.6,5.2);
			\node at (2,5.3)  { $d_{11}$};
			\draw[ thick, black](2.4, 5.2) -- (3,5.2);
			\draw[ thick, black](3, 5.1) -- (3,5.4);
			\draw[ thick, black](3, 5.2) -- (3.6,5.2);
			\node at (4,5.3)  { $d_{11}$};
			\draw[ thick, black](4.4, 5.2) -- (5,5.2);
			\draw[ thick, black](5, 5.1) -- (5,5.4);
			
			\draw[ thick, black](6.1, 5) -- (6.4,5);
			\draw[ thick, black](6.3, 5) -- (6.3,4.7);
			\node at (6.3,4.6)  { $d_{12}$};
			\draw[ thick, black](6.1, 4) -- (6.4,4);
			\draw[ thick, black](6.3, 4.4) -- (6.3,4);
			\draw[ thick, black](6.3, 4) -- (6.3,2.6);
			\node at (6.3,2.5)  { $H$};
			\draw[ thick, black](6.3, 2.3) -- (6.3,1);
			\draw[ thick, black](6.1, 1) -- (6.4,1);
			\draw[ thick, black](6.3, 1) -- (6.3,0.7);
			\node at (6.3,0.6)  { $d_{12}$};
			\draw[ thick, black](6.3, 0.4) -- (6.3,0);
			\draw[ thick, black](6.1, 0) -- (6.4,0);

			\filldraw[fill=blue, draw=black] (13,1) rectangle (18,5);
			\draw[step=1,black, thick] (13,1) grid (18,5);
			\draw[step=1,black, thick] (14,0) grid (17,6);
			
			\draw[ultra thick, red](14, 3) -- (17,3);
			\node at (15.5,2.5)  { \huge $E_{i}$} ;
			\node at (3,-0.5)  { \huge $\omega_{i}^+$} ;
			\draw[ thick, black](14,-0.1) -- (14,-0.4);
			\draw[ thick, black](13,-0.1) -- (13,-0.4);
			\draw[ thick, black](17,-0.1) -- (17,-0.4);
			\draw[ thick, black](18,-0.1) -- (18,-0.4);
			\draw[ thick, black](13, -0.2) -- (13.2,-0.2);
			\node at (13.5,-0.2)  { $d_{21}$};
			\draw[ thick, black](13.7, -0.2) -- (14,-0.2);
			\draw[ thick, black](14, -0.2) -- (15.2,-0.2);
			\draw[ thick, black](15.7, -0.2) -- (17,-0.2);
			\node at (15.5,-0.2)  { $H$};
			\draw[ thick, black](17, -0.2) -- (17.2,-0.2);
			\draw[ thick, black](17.7, -0.2) -- (18,-0.2);
			\node at (17.5,-0.2)  { $d_{21}$};
			
			\draw[ thick, black](18.1,1) -- (18.4,1);
			\draw[ thick, black](18.1,3) -- (18.4,3);
			\draw[ thick, black](18.1,5) -- (18.4,5);
			\draw[ thick, black](18.3,5) -- (18.3,4.1);
			\draw[ thick, black](18.3,3.7) -- (18.3,3);
			\node at (18.3,3.9)  { $d_{22}$};
			\draw[ thick, black](18.3,3) -- (18.3,2.1);
			\draw[ thick, black](18.3,1.7) -- (18.3,1);
			\node at (18.3,1.9)  { $d_{22}$};
			\end{tikzpicture}
		}
		\caption{ Illustration of an oversampled neighborhood associated with the coarse edge $E_i$.}
		\label{fig:oversampgrid}
	\end{figure}

 {\bf Step 2:}  For each coarse face $E_i$, we select an oversampled
 coarse neighborhood $\omega^+_i$ (see Figure \ref{fig:oversampgrid} for illustration). We repeat the selection for
 coarse faces $E_1, E_2, \cdots E_I$. Then we obtain oversampled neighborhoods $\omega_1^+$, $\omega_2^+ $,
 $\cdots$, $\omega_I^+$ $\subseteq D$.
  The index $\{1,2,...,I\}$ can be chosen such that $\displaystyle\bigcup_{i=1,\cdots,I}\omega_i$ (not $\omega_i^+$)
 form a non-overlapping partition of $D$.

{\bf Step 3:}   For each $\omega_i^+$, we solve for $\phi_i^+\in \widehat{V}^{i,+}_{\text{snap}}$ such that
 \begin{equation}\label{eq:residual equation}
\int_{\omega_i^+}\kappa^{-1}\phi_i^+\cdot v=R_{\omega_i^+}(v)\quad\quad\quad \forall v\in \widehat{V}^{i,+}_{\text{snap}}
 \end{equation}
  We note that we need to project the $v_{H}^l$ and $p_{H}^l$ into the fine-grid spaces $V_h$ and $Q_h$, respectively, in order to compute $R_{\omega_i^+}(v)$.

{\bf Step 4:}  We take the restriction of $\phi_i^+\cdot n_i^+$ on the coarse face $E_i$ and normalize it, and we denote it
 by $\lambda_{i}$.

 {\bf Step 5:} The online basis $\chi_i$ for the coarse neighborhood $\omega_i$ can be constructed by solving
 \begin{equation}\label{eq:compute_onlinevelocity}
 \begin{split}
 \kappa^{-1}\chi_i+\nabla p_i&=0 \quad \text{in } \omega_i,\\
 \text{div}(\chi_i)&=\alpha_i\quad \text{in } \omega_i,\\
 \chi_i\cdot n_i&=\lambda_{i}\quad \text{on } E_i,\\
  \chi_i\cdot n_i&=0\quad \text{on }\partial\omega_i,
 \end{split}
 \end{equation}
 where $\alpha_i$ is chosen to satisfy the condition $\int_K\alpha_i=\int_{\partial K}v_i\cdot n_i$
 for every $K\subseteq\omega_i$,
 $n_i$ is a fixed unit-normal vector for the coarse face $E_i$.
 Those $\{\chi_i\}_{i=1}^I$ are the new online basis functions.  We update the velocity
 space by letting $V_H^{l+1}=V_H^{l}\oplus \text{span}\{\chi_1,\chi_2,\cdots,\chi_I\}$.

 After Step 5, we repeat from  Step 1 until we have certain number of basis functions or some predefined error indicator is small.

Next we present the construction of a local online basis function for
coarse face $E_i$ in matrix formulation at level $l$.
We define
\begin{equation*}
R^l=[\Psi_1^l,\Psi_2^l,\cdots,\Psi_{Nl}^l]
\end{equation*}
where $\{\Psi_i^l\}_{i=1}^{Nl}$ are the functions in space $V_H^{l}$.
We also define a $J_i^+\times(J_i^+-1) $ matrix
\[
P_i^+
=
\begin{bmatrix}
1 & 0 & 0 &0 & \dots  & 0 \\
-1 & 1 & 0 &0 & \dots  & 0 \\
\vdots & \vdots & \vdots & \ddots & \vdots& \vdots \\
0 & 0 & -1&1 & 0  &\vdots\\
0 & 0 & \dots&-1 & 1  &0\\
0 & 0 & 0 &\dots & -1  &1
\end{bmatrix}
\]
This matrix is to find the divergence free subspace $\widehat{V}^{i,+}_{\text{snap}}$ of ${V}^{i,+}_{\text{snap}}$.
We remark that the definition of $P_i^+$ is not unique.

Then the online basis construction algorithm in matrix formulation is as follows:

 {\bf Step 1:}  Solve  the coarse system
\begin{equation}\label{eq:coarse_sys}
\left[\begin{array}{cc}
(R^l)^TA_hR^l & (R^l)^TB_h^{T}G_H \\
G_H^TB_hR^l & 0\\

\end{array}\right]\left[\begin{array}{cc}v_{Hr}^l\\p_{Hr}^l\end{array}\right] =\left[\begin{array}{cc}0\\G_H^TF_h\end{array}\right]
\end{equation}
where $v_{Hr}^l$ and $p_{Hr}^l$ are vectors of coefficients in the expansions
of the solutions $v_H$, $p_H$ in the space $V_H^l$ and $Q_H$.
Note $v_H^l$ is defined in space $V_H^l$, we can project it into space
$v_h$ by computing $v_{Hr}^{h,l}=R^lv_{Hr}^l$.
We also need to compute the projection of $p_H^l$ in space
$P_h$ by using $p_{Hr}^{h,l}=G_Hp_{Hr}^l$.

{\bf Step 2:}  Compute the all local residuals $r_i^+$ by $r_i^+=\big(F_h-A_hv_{Hr}^{h,l}-B_h^Tp_{Hr}^{h,l}\big)\vert_{\omega_i^+}$.

{\bf Step 3:}  Solve the residual Equation (\ref{eq:residual equation}), that is to solve
\begin{equation*}
\bigg((P_i^+)^TA_{\text{snap}}^{i,+}P_i^+\bigg)v_{hr}^{i,+}=(P_i^+)^T(R_{\text{snap}}^{i,+})^Tr_i^+
\end{equation*}
where $A_{\text{snap}}^{i,+}=(R_{\text{snap}}^{i,+})^TA_h^{i,+}R_{\text{snap}}^{i,+}$, which can be precomputed before online iterations and saved. The dimension of this matrix is $J_i^+\times J_i^+ $.
Note that $v_{hr}^{i,+}\cdot n_i^+$ lies on $E_i^+$, take the restriction of it on $E_i$ and normalize it, then we get
 $\lambda_i$ in the system (\ref{eq:compute_onlinevelocity}).

{\bf Step 4:}  We can compute the online basis $\chi_i$ by $\chi_i=R_{\text{snap}}^{i}\lambda_{i}$
if $R_{\text{snap}}^{i}$ is precomputed. If not, one can use the lowest-order Raviart-Thomas finite element method
to solve the system (\ref{eq:compute_onlinevelocity}) to get the online basis for $\omega_i$.

\section{Numerical examples}\label{numerical-results}
In this section, we present several representative examples to show the performance of our oversampled online method.
In all simulations reported below, the computational domain is $D = (0,1)^d$.
 The fine-grid solution is used as the reference solution in all numerical examples.
 We use the preprocessing steps introduced in \cite{mixed_dd1988} for coarse-grid problem to remove the effects of singular source used in practical
 flow simulations.

\begin{figure}[H]
	\centering
	\subfigure[model 1: $\kappa_{1}$]{
		\includegraphics[width=3.4in]{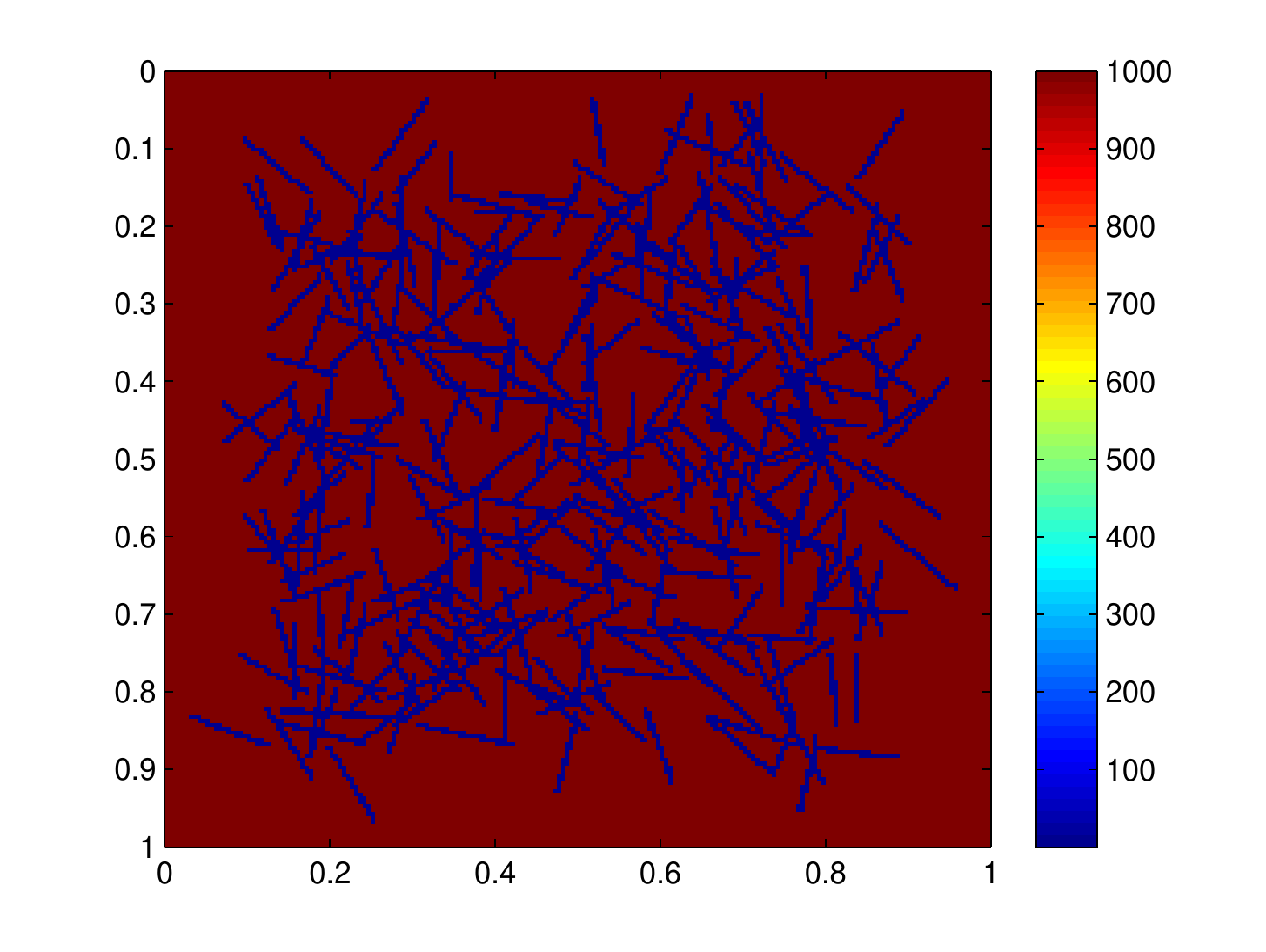}}
		\subfigure[model 2: $\kappa_{2}$ in $\log_{10}$ scale]{
			\includegraphics[width=3.4in]{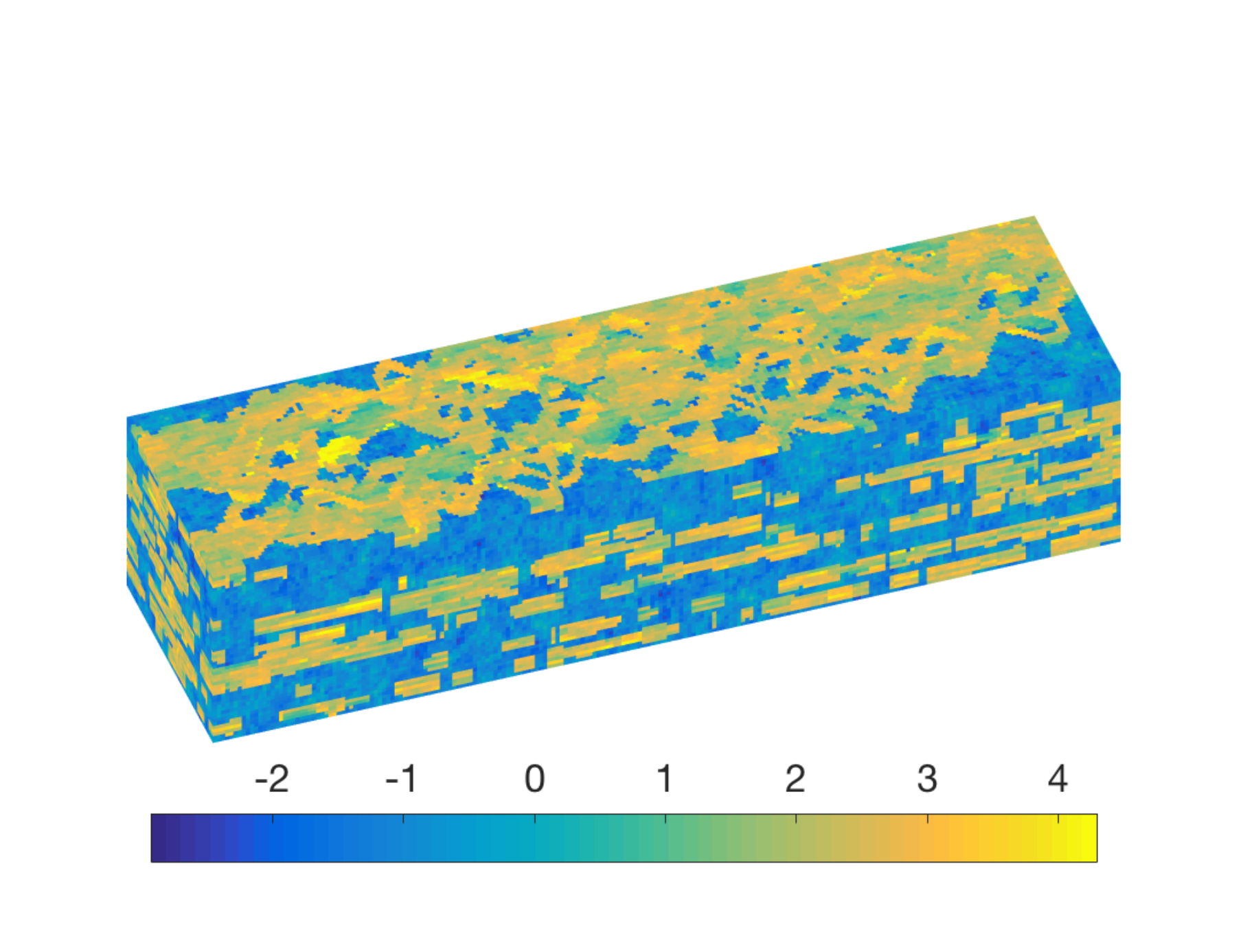}}
				\subfigure[model 3: $\kappa_{3}$ ]{
			\includegraphics[width=3.4in]{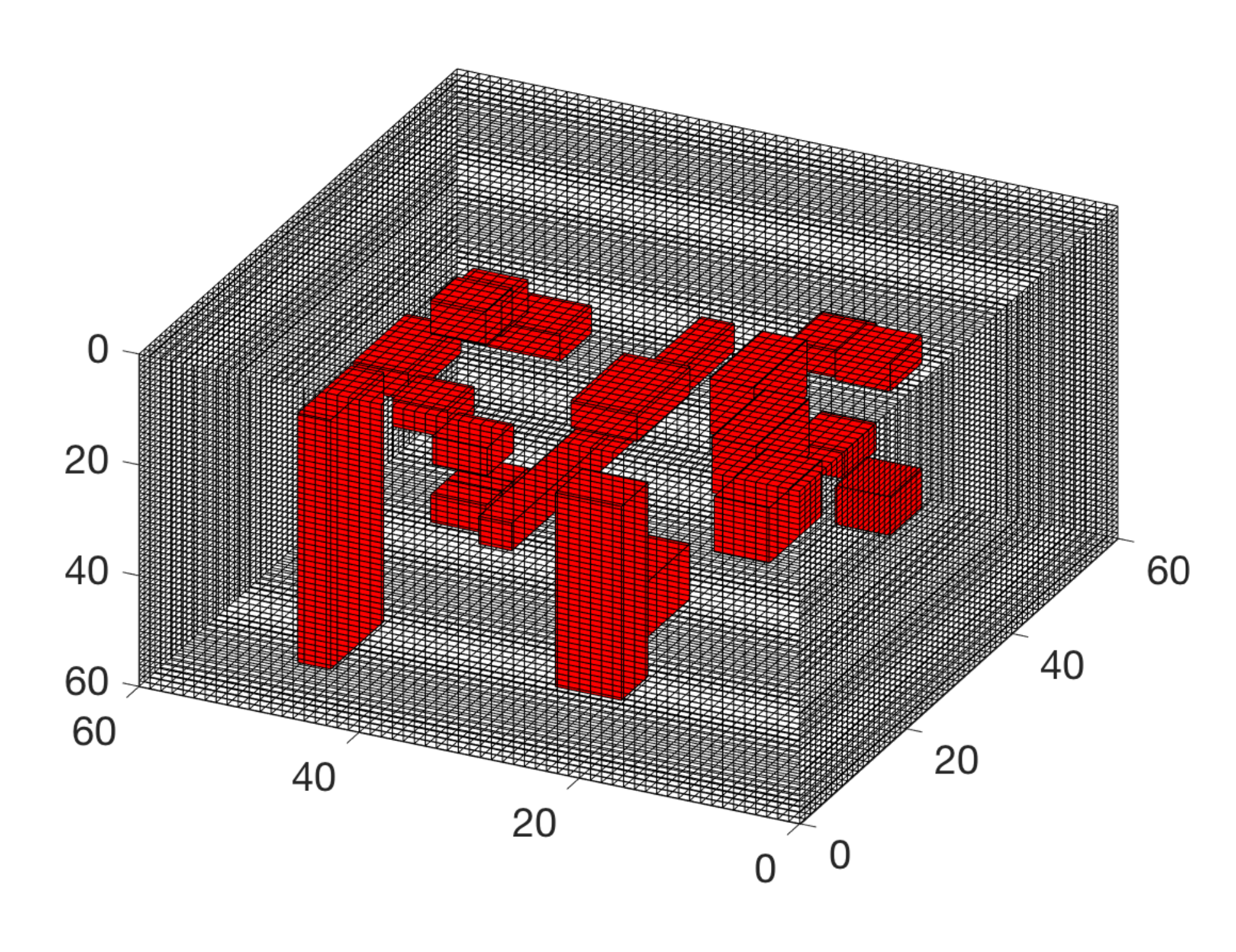}}
	\caption{Three permeability fields in the numerical examples.}
	\label{fig:model} 
\end{figure}

 We consider three models with permeability $\kappa$ depicted in Figure \ref{fig:model}.  We note that for model 1, $\kappa=0.1$ in the blue region and $10^3$ in the red region.
As it is shown, model 3 contains high contrast, long channels, and isolated inclusions, we set the value in the blank region equals 1 and the red region equals $k_0$. Model (b) contains the last 50 layers of  the SPE10 model \cite{Aarnes2005257}.  The fine grid mesh for this model is
$220 \times 60 \times 50$. The SPE 10 model is  used as a  benchmark  to test different upscaling techniques and multiscale methods,  and is therefore a good test case for our methodology.
For all simulations presented below, we divide the computational domain $D$ into $N_x\times N_y$
($N_x\times N_y\times N_z $ in 3D) coarse blocks, and in each coarse block, we further generate
a uniform $n\times n$ ($n\times n\times n$) fine scale square elements.
We define the following weighted velocity error to quantify the accuracy of the online multiscale solution
\begin{equation*}
 e_{\boldsymbol{v}}:=\frac{\|\boldsymbol{v}_{ms}-\boldsymbol{v}_f\|_{\kappa,D}}{\|\boldsymbol{v}_f\|_{\kappa,D}}
\end{equation*}
where $\|\boldsymbol{v}\|_{\kappa,D}^2=\int_{D}\kappa^{-1}\boldsymbol{v}^2$.

We test our method for both elliptic problem and two-phase flow simulation. First, we show the results for the elliptic problem in Section \ref{num exam: elliptic}.

\begin{figure}[H]
	\centering
	\subfigure{
		\includegraphics[width=3.5in]{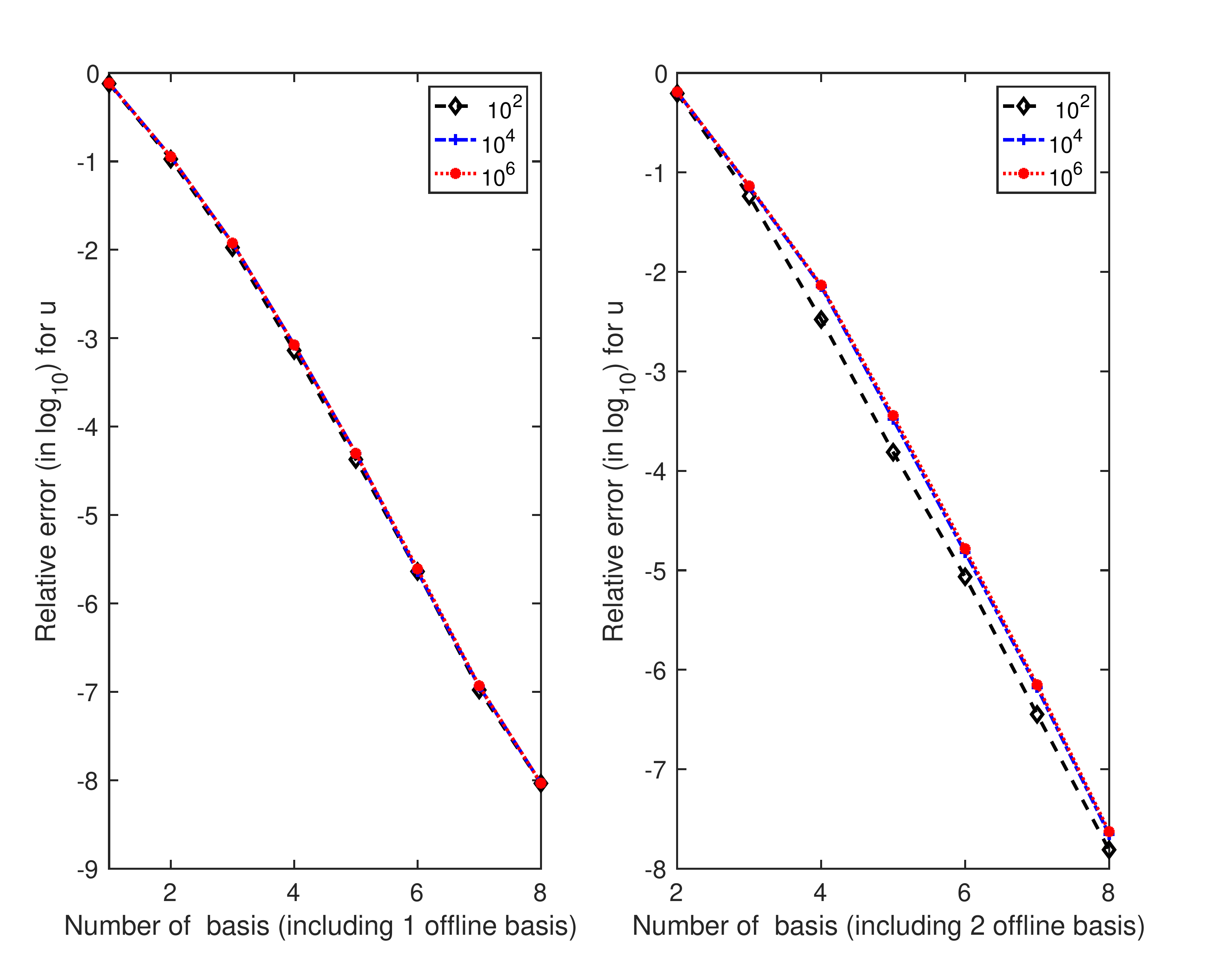}}
	\caption{Relative weighted velocity error with $k_0=10^2, 10^4, 10^6$,
		model 3.}
	\label{fig:contrast1} 
\end{figure}

\begin{figure}[H]
	\centering
		\includegraphics[width=3.5in]{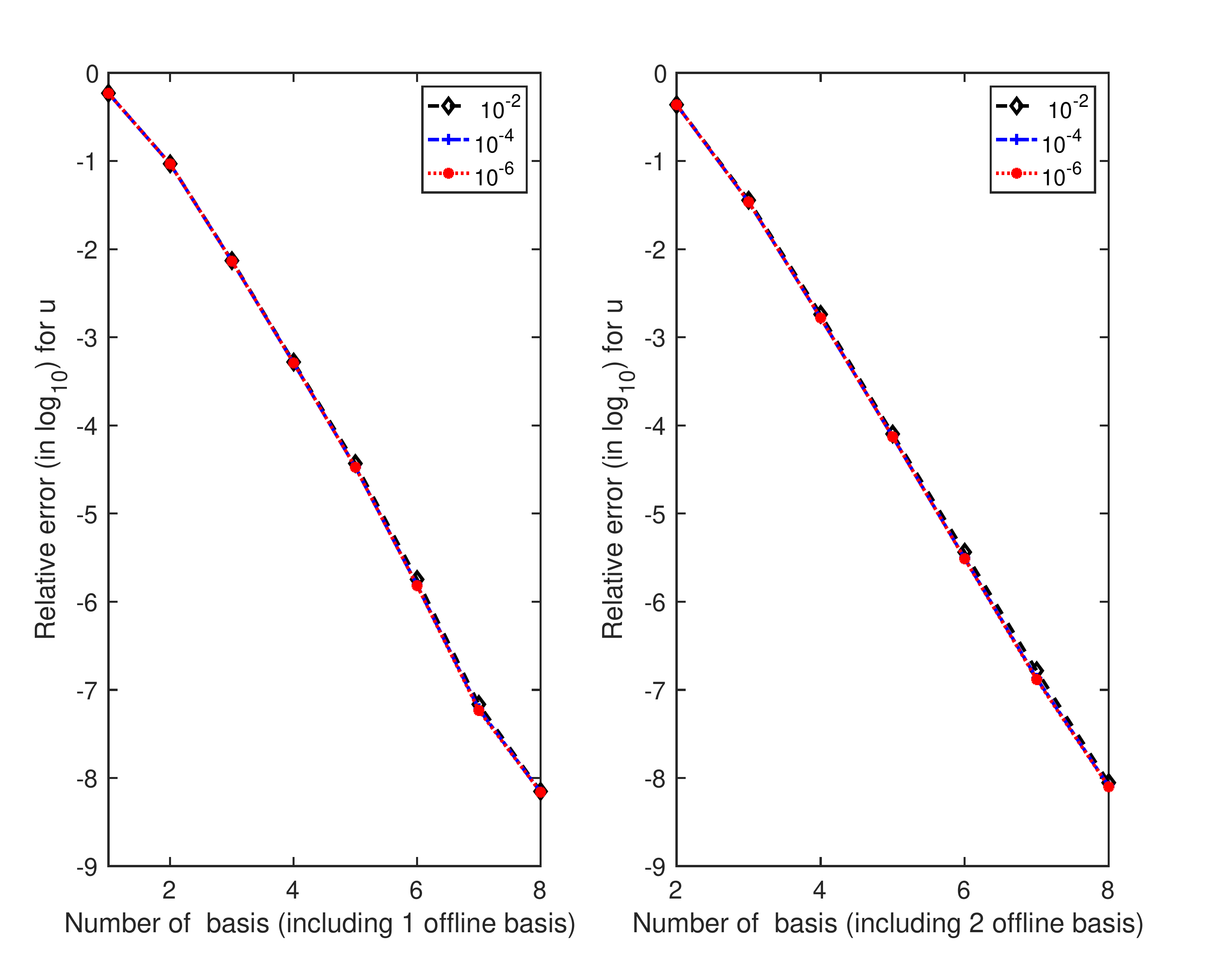}
	\caption{Relative weighted velocity error with $k_0=10^{-2}, 10^{-4}, 10^{-6}$, model 3.}
	\label{fig:contrast2} 
\end{figure}

\subsection{Oversampled online method for elliptic problems}\label{num exam: elliptic}

In this section, we first compare the performance of  online enrichment with and without oversampling, and show that oversampling helps achieve faster convergence. We also discuss the influence of oversampled domain size on convergence.  We compare the efficiency of the online method for the elliptic problem by using  different contrast orders for both low and high contrast orders and show that convergence is almost independent of contrast orders.

We use $\left(\begin{array}{cc}
d_{11} & d_{12} \\
d_{21}  & d_{22} \\
 \end{array}
 \right)$ to define  the local domain (see Figure~\ref{fig:oversampgrid}) for the computation of the online basis functions for 2D computation. In total, three cases  which are given below are considered.  Here as defined earlier,  $n$ is the number of fine elements in a coarse block for each direction. The local domain in Case 1 is exactly the neighborhood of  an coarse edge, i.e., no oversampling. Case 2 and Case 3 are oversampling domains. In Case 2, one layer of fine cell is added to the coarse neighborhood in the direction
of the edge. Therefore, the oversampling domain is larger than the coarse neighborhood. In Case 3, the
layers of cells on both sides of the coarse edge are reduced to about a half of n, while one layer of fine
cells is added to the coarse neighborhood in the direction of this edge. Thus, the oversampling domain
in this case is smaller than the coarse neighborhood.\\
Case 1: No oversampling: $\left(\begin{array}{cc}
$n$ & 0 \\
0 & $n$\\
\end{array}
\right)\\$
Case 2: oversampling case a: $\left(\begin{array}{cc}
$n$ & 1 \\
1& $n$\\
\end{array}
\right)\\$\\
Case 3: oversampling case b: $\left(\begin{array}{cc}
$ [n/2] $ & 1 \\
1 & $ [n/2] $\\
\end{array}
\right)\\$\\

In each table of Tables \ref{model-a-err-1}-\ref{model-c-err-2}, we compare the velocity errors of the above 3 cases for each test model. Tables \ref{model-a-err-1}-\ref{model-a-err-2} are results for model 1, with $N_x=10$,$N_y=10$, $n=20$. In Table \ref{model-a-err-1}, we start with 1 initial
offline basis function in each coarse neighborhood $\omega_i$, then we iteratively add online basis functions up to 7. The first column gives the number of basis functions represented by $N_b$ in each coarse neighborhood and the dimension of the coarse system (\ref{eq:coarse_sys}). The second column is the error decay history for Case 1, which as we mentioned earlier is the non-oversampling case. The third and the fourth columns are for Case 2 and Case 3, respectively.
It can be found that with more online basis functions  added, the multiscale solutions from  the three cases all
converge to the reference solution. However,
one can observe that the online solution with oversampling converges much faster
than the no-oversampling case. For example, if 7 basis functions are used for the no-oversampling case, the error is around $10^{-3}$, while for the Case 2 or 3
one can obtain similar accuracy solution with only 4 basis functions.
We remark that adding more fine-grids in the oversampling
domain one can yield faster error decay.
Though the size of oversampled domain for Case 3 is much smaller for
Case 2, however, its corresponding convergence rate is comparable with
Case 2 since
one layer of fine grid cell is added on each side in the direction of coarse
edge (face).
Convergence behavior for each case with two initial offline basis functions is shown in Table \ref{model-a-err-2}.  By comparing corresponding columns for the three cases, one can see that by using one more initial basis function, faster convergence is achieved. For example, by comparing the columns for Case 2 in Table \ref{model-a-err-1} and Table \ref{model-a-err-2}, the rate of convergence is faster in Table \ref{model-a-err-2}. The error convergence history for both model 2 and  model 3
are reported in Tables \ref{model-b-err-1}-\ref{model-c-err-2}.
We can draw similar conclusion as for model 1.

Next we study the performance of different $k_0$ for model 3.
Figure \ref{fig:contrast1} shows the  errors against the number of
basis functions for $k_0=10^2,10^4,10^6$. The left figure is the
case of using 1 initial offline basis while the right figure is for
using 2 initial basis functions.
We can see that for both cases the error convergence history is almost the
same for different $k_0$. We also test the lower contrast case, i.e., $k_0=10^{-2},10^{-4},10^{-6}$.
Corresponding results are plotted in Figure \ref{fig:contrast2},
we found similar phenomenon, that is, the difference in the contrast has almost no effects on the relative errors.
We can conclude that the oversampled online enrichment method is robust in the sense that its convergence rate is almost 	independent of the media's contrast
order.

\begin{table}[H]
	\centering \begin{tabular}{|c|c|c|c|c|c|c}
		\hline
		$N_b$(Dim)& Case 1   & Case 2&Case 3  \tabularnewline	\hline
		1(280) &1.06e+0&1.06e+0&1.06e+0 \tabularnewline\hline
		2(460)&1.85e-01&7.38e-02&1.13e-01  \tabularnewline\hline
		3(640) &2.94e-02 &9.67e-03& 1.31e-02\tabularnewline\hline
		4(820) &1.18e-02&1.38e-03 & 1.76e-03 \tabularnewline\hline
		5(1000) &1.05e-02&2.40e-04&2.30e-04  \tabularnewline\hline
		6(1180) &2.48e-03&1.25e-06 &7.35e-06\tabularnewline\hline
		7(1360)&1.22e-03 &2.40e-08 &4.51e-08   \tabularnewline\hline
	\end{tabular}
	\caption{1 initial offline basis, model 1, $N_x=10$, $N_y=10$, $n=20$.}
	\label{model-a-err-1}
\end{table}

\begin{table}[H]
	\centering \begin{tabular}{|c|c|c|c|c|c|c}
		\hline
		$N_b$(Dim)& Case 1   & Case 2&Case 3  \tabularnewline	\hline
		2(460) &9.51e-01&9.51e-01& 9.51e-01  \tabularnewline\hline
		3(640) &1.38e-01  &5.38e-02 & 6.19e-02   \tabularnewline\hline
		4(820)&1.38e-02&1.26e-03 & 2.21e-03 \tabularnewline\hline
		5(1000) &4.05e-03&8.39e-05 &3.54e-04 \tabularnewline\hline
		6(1180) &2.58e-03 &3.52e-06& 1.28e-05  \tabularnewline\hline
		7(1360) &6.41e-04 &4.73e-08 & 1.08e-07  \tabularnewline\hline
	\end{tabular}
	\caption{2 initial offline basis, model 1, $N_x=10$, $N_y=10$, $n=20$ .}
	\label{model-a-err-2}
\end{table}


\begin{table}[H]
	\centering \begin{tabular}{|c|c|c|c|c|c|c}
		\hline
		$N_b$(Dim)& Case 1   & Case 2&Case 3  \tabularnewline	\hline
		1(2368) &6.54e-01&6.54e-01&6.54e-01 \tabularnewline\hline
		2(4076) &1.03e-01&5.22e-02&7.07e-02   \tabularnewline\hline
		3(5784) &3.23e-02&8.13e-03&1.11e-02 \tabularnewline\hline
		4(7492) &1.73e-02&2.58e-03&3.02e-03  \tabularnewline\hline
		5(9200) &1.35e-02&4.10e-04& 5.78e-04 \tabularnewline\hline
		6(10908) &1.00e-02&2.38e-05&3.26e-05 \tabularnewline\hline
		7(12616) &7.92e-03&1.17e-06& 2.74e-06   \tabularnewline\hline
		8(14324) &6.07e-03&4.09e-08&1.32e-07  \tabularnewline\hline
	\end{tabular}
	\caption{1 initial offline basis, model 2, $N_x=6$, $N_y=22$, $N_z=5$, $n=10$.}
	\label{model-b-err-1}
\end{table}

\begin{table}[H]
	\centering \begin{tabular}{|c|c|c|c|c|c|c}
		\hline
		$N_b$(Dim)& Case 1   & Case 2&Case 3  \tabularnewline	\hline
		2(4076) &6.34e-01&6.34e-01& 6.34e-01 \tabularnewline\hline
		3(5784) &8.91e-02&4.39e-02&5.72e-02 \tabularnewline\hline
		4(7492) &2.23e-02&3.20e-03&5.09e-03  \tabularnewline\hline
		5(9200) &5.52e-03&1.72e-04&3.10e-04 \tabularnewline\hline
		6(10908) &2.53e-03&1.21e-05&2.80e-05 \tabularnewline\hline
		7(12616) &1.85e-03&6.21e-07&1.85e-06  \tabularnewline\hline
		8(14324) &1.21e-03&2.87e-08&8.77e-08  \tabularnewline\hline
	\end{tabular}
	\caption{2 initial offline basis, model 2, $N_x=6$, $N_y=22$, $N_z=5$, $n=10$.}
	\label{model-b-err-2}
\end{table}

\begin{table}[H]
	\centering \begin{tabular}{|c|c|c|c|c|c|c}
		\hline
		$N_b$(Dim)& Case 1   & Case 2&Case 3  \tabularnewline	\hline
		1(756) &7.69e-01&7.69e-01&7.69e-01 \tabularnewline\hline
		2(1296) &1.12e-01&6.99e-02&1.13e-01  \tabularnewline\hline
		3(1836) &2.88e-02&5.37e-03&1.19e-02 \tabularnewline\hline
		4(2376) &1.04e-02&2.94e-04&8.40e-04  \tabularnewline\hline
		5(2916) &3.74e-03&1.67e-05&4.98e-05 \tabularnewline\hline
		6(3456) &1.23e-03&7.13e-07&2.47e-06 \tabularnewline\hline
		7(3996) &4.04e-04&2.96e-08&1.30e-07   \tabularnewline\hline
		8(4536) &2.69e-04&9.36e-10&4.50e-09  \tabularnewline\hline
	\end{tabular}
	\caption{1 initial offline basis, model 3, $k_0=10^4$, $N_x=6$, $N_y=6$, $N_z=6$, $n=10$.}
	\label{model-c-err-1}
\end{table}

\begin{table}[H]
	\centering \begin{tabular}{|c|c|c|c|c|c|c}
		\hline
		basis(Dim)& Case 1   & Case 2&Case 3  \tabularnewline	\hline
		2(1296) &6.30e-01&6.30e-01& 6.30e-01 \tabularnewline\hline
		3(1836) &1.03e-01&3.91e-02&5.77e-02 \tabularnewline\hline
		4(2376) & 2.59e-02&1.91e-03&3.40e-03  \tabularnewline\hline
		5(2916) &7.01e-03&7.76e-05&1.55e-04 \tabularnewline\hline
		6(3456) &1.67e-03&3.99e-06&8.34e-06 \tabularnewline\hline
		7(3996) &4.13e-04&1.61e-07&3.51e-07   \tabularnewline\hline
		8(4536) &8.18e-05&5.11e-09&1.24e-08  \tabularnewline\hline
	\end{tabular}
	\caption{2 initial offline basis, model 3, $k_0=10^4$, $N_x=6$, $N_y=6$, $N_z=6$, $n=10$.}
	\label{model-c-err-2}
\end{table}

\subsection{A two phase flow and transport problem}
In this section, we test our method by solving an incompressible two-phase flow and transport model problem, which is
 used to simulate porous media flows \cite{efendiev2016online,SpeJ}.  The mixed formulation for the elliptic problem
 is particularly suitable for two-phase flow problems. As the velocity field is mass conservative, which is an important feature
 to get accurate saturation solution.
 In particular, we consider two-phase flow in a reservoir domain (denoted by $\Omega$)
with the assumption that the fluid displacement is only driven by viscous effects, neglecting compressibility
and gravity for simplicity. The two phases we  consider are water and oil, which are assumed to be immiscible. Next, we summarize the differential equations for the two-phase flow. By the Darcy's law, we have the following equation
for each phase
\begin{equation} \label{darcy}
{ \bf v}_l=-\frac{k_{rl}(s_l)}{\mu_l} { K} \nabla {p}
\end{equation}
where ${  \bf v}_l$ is the phase velocity, $ {K}$ is the permeability
tensor, $k_{rl}$ is the relative permeability to phase $l$ ($l=o, w$), $s_l$ is
 saturation, and ${p}$ is pressure.
Throughout the paper, we use a
single set of relative permeability.

By the mass conservation law, the following equation for each phases is obtained:
\begin{equation} \label{mass}
\phi\frac{\partial s_l}{\partial t} + \nabla \cdot {\bf v}_l = q_l.
\end{equation}

To close the system,  the property $s_w+s_o=1$ is used. Then, we can  write the whole system depending on pressure and saturation as follows (we use $s$ instead of $s_w$ for simplicity):
\begin{eqnarray}
\nabla\cdot {\bf v} &=& q_w + q_o \quad \textrm{in} \quad \Omega \label{pressEq}\\
\phi \frac{\partial s}{\partial t}+ \nabla \cdot({ f_w(s) {\bf v}}) &=& \frac{q_w}{\rho_w} \quad \textrm{in} \quad \Omega \label{satEq}\\
{\bf v}\cdot n & = & 0 \quad \textrm{on} \quad \partial {\Omega} \quad \textrm{(no flow at boundary)}\\
s(t=0) & = & s_{0} \quad \textrm{in} \quad \Omega \quad \textrm{(initial known saturation)} \label{pressEq4}
\end{eqnarray}
where $\phi$ is the porosity, $\lambda$ is the total mobility defined as
\begin{equation} \label{mobility}
\lambda(s) =\lambda_w(s)+\lambda_o(s)=\frac{k_{rw}(s)}{\mu_w}+\frac{k_{ro}(s)}{\mu_o}
\end{equation}
 $f_w(s)$ is the flux function,
\begin{equation} \label{fraction}
f_w(s) =\frac{\lambda_{w}(s)}{\lambda(s)}=\frac{k_{rw}(s)}{k_{rw}(s)+\frac{\mu_w}{\mu_o}k_{ro}(s)}
\end{equation}
and ${\bf v} = { \bf v}_w + { \bf v}_o=-\lambda(s) { K}  \nabla {p}$ is the total flux. Moreover, $q_w$ and $q_o$ are volumetric source terms
for water and oil.

Here, sequential formulation is employed to solve the coupled system (\ref{pressEq})-(\ref{pressEq4}), that is at each time step   the pressure equation is solved first to get the velocity, then  the updated velocity is used to solve the saturation equation. The pressure equation is solved by using the offline basis functions together with the online basis functions  computed at the initial time step,  and the saturation equation is solved by finite volume method.

 The initial saturation is taken to be zero. The velocity in Equation (\ref{satEq}) is the fine grid velocity which is obtained by projecting the multiscale velocity field  onto the fine grid.
  We use a five-spot model, i.e., we  inject the water in the four corners of the model, and  the sink is located at the center of the model.
 We test the method for model 1 and 2 with different
 combination of the number of  offline and online basis functions.
  Time step size is 50, and the end of simulation time is 5000. We define the relative $L^2 $ saturation error at time step $i$ as\begin{equation*}
e_s(i):=\frac{\|s_{ms}(i)-s_{ref}(i)\|_{L^2,\Omega}}{\|s_{ref}(i)\|_{L^2,\Omega}}.
\end{equation*}
The average of $e_s$ is simply the arithmetic average of $e_s(i), i=1,\cdots, 5000$.

\begin{figure}[H]
	\centering
	\subfigure{
		\includegraphics[scale=0.30]{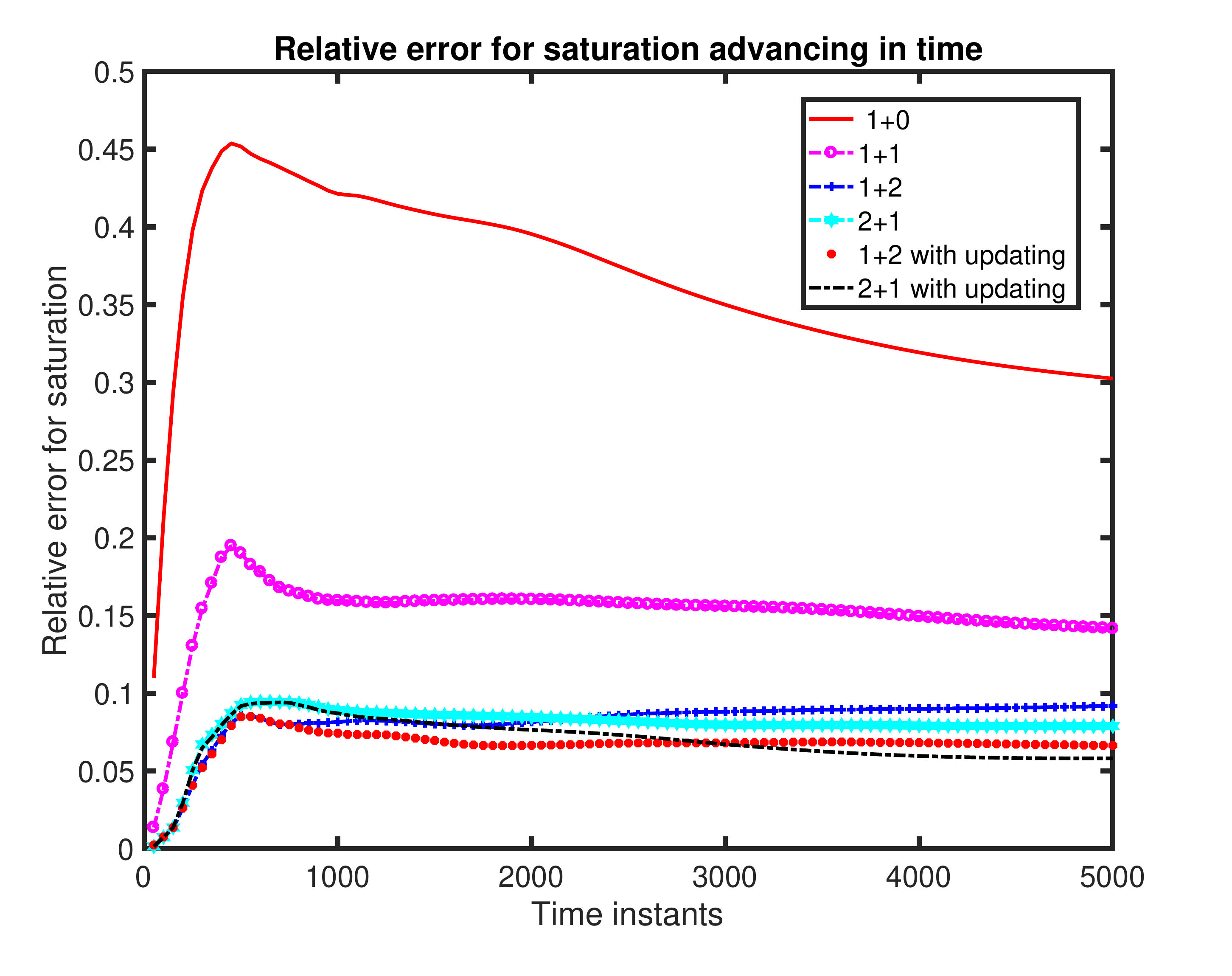}}	
		
	\caption{Saturation error for model 1.}	
	\label{fig:Sat_model_frac}
\end{figure}

\begin{figure}[H]
	\centering
	\subfigure[]{
		\includegraphics[scale=0.30]{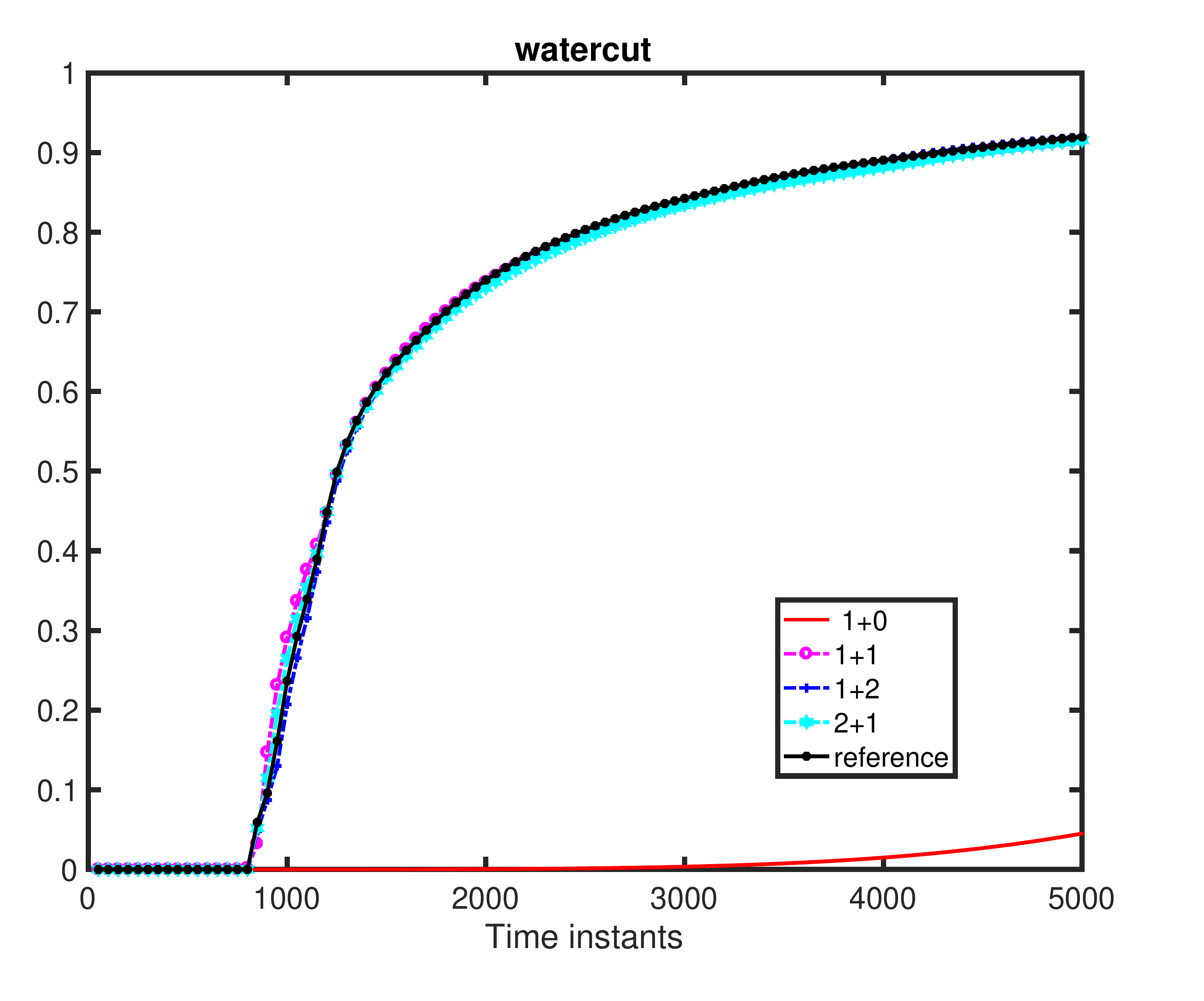}}	
	\subfigure[]{
		\includegraphics[scale=0.30]{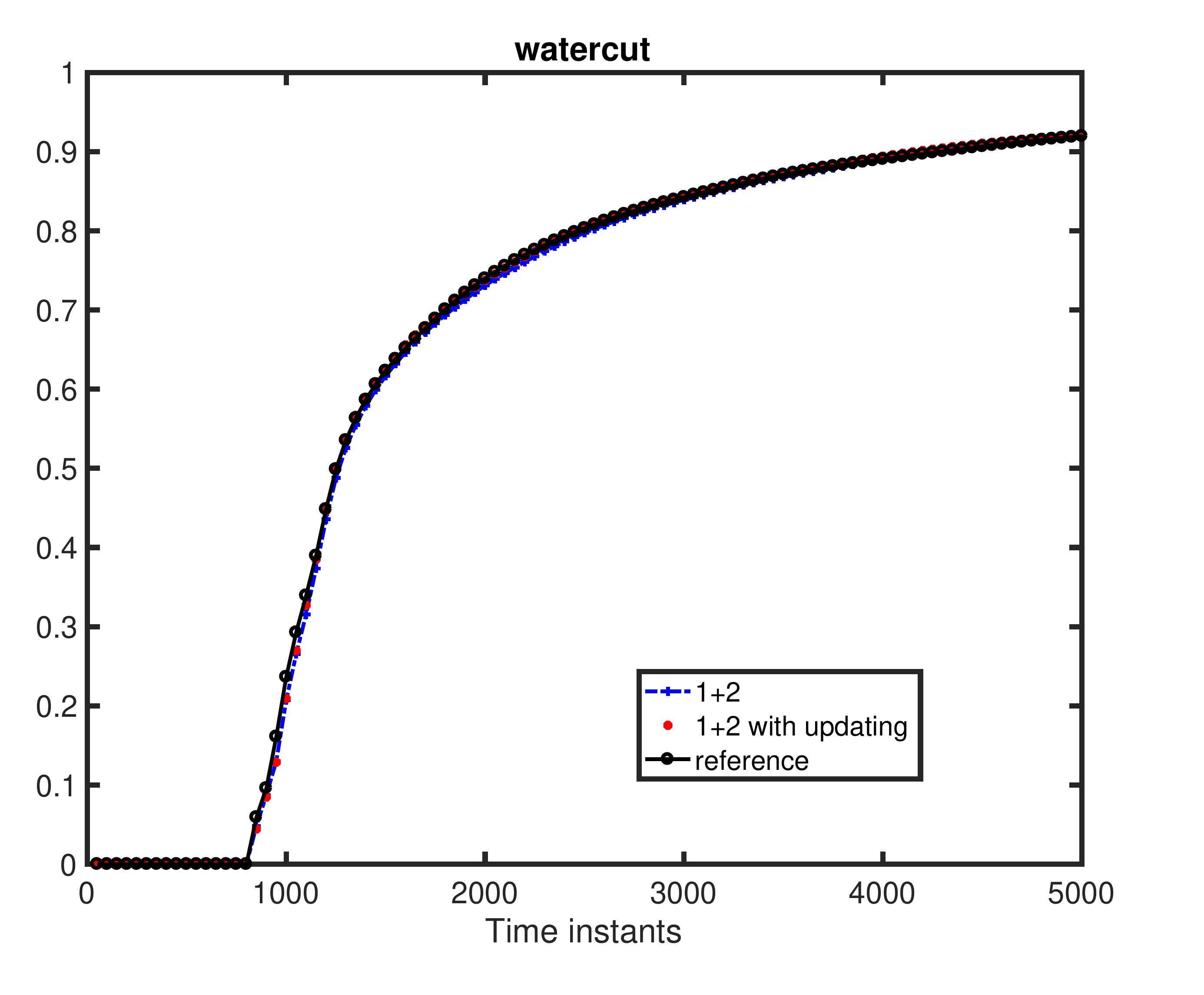}}	
	\subfigure[]{
		\includegraphics[scale=0.30]{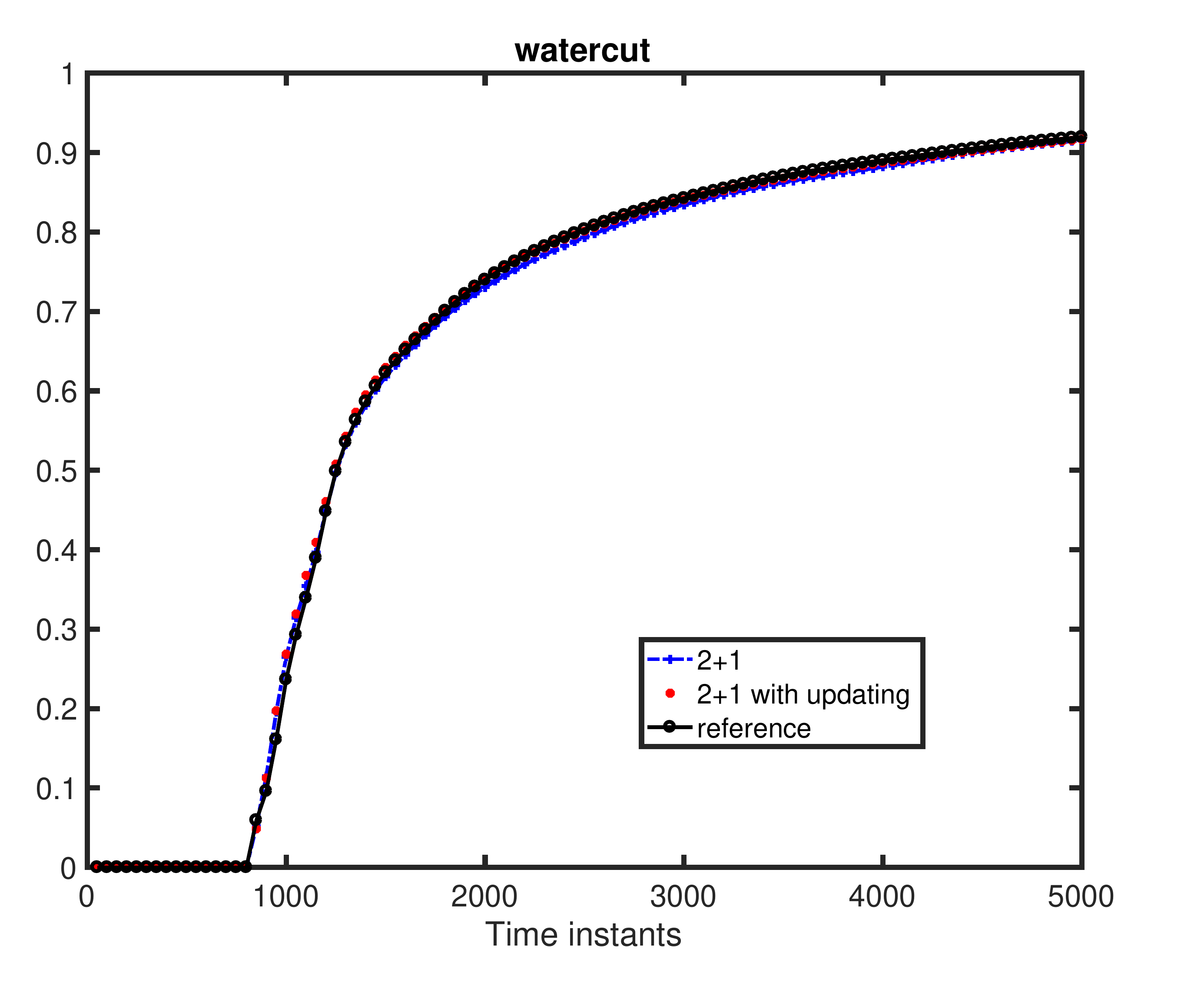}}

		
	\caption{Water-cut for model 1.}	
	\label{fig:watercut_frac}
\end{figure}

Figure \ref{fig:Sat_model_frac} shows the relative saturation errors of using different online basis functions and different number of offline basis
function against the time instants for model 1.
 We note that $x + y$ in the legend of the figure means using
 $x$ offline basis functions and $y$ online basis functions.
 Updating in the legend means that we update the velocity basis functions
at every  other five time steps
 by solving local problems on the coarse neighborhoods
 with the initial flux on the coarse face as the boundary condition.
 The updated basis may contain features of the updated permeability field,
 and thus yields more accurate saturation.
 Otherwise we use the online basis functions obtained with initial model.
  We observe the errors are greater than $30\%$ for almost all time instants if only 1 offline basis is used (shown by the red line). However, if 1 or 2 online basis functions are added (the purple and blue lines), the error can be reduced tremendously.
  Updating the basis can further decrease the saturation errors.
  We also present the water-cut comparison (water flux fractional function
  $f_w(s))$ in Figure \ref{fig:watercut_frac}.
  We can see that the water cut lines for all the cases that have at least
  1 online basis are very close to the reference line especially
  when time-instant is larger than 2000. However, if we do not update
  the basis, there is noticeable error between multiscale water-cut and
  the reference water-cut. Updating the basis function can make the multiscale
  water-cut line and reference water-cut line almost identical.

\begin{figure}[H]
	\centering
	\subfigure[Fine-scale solution ]{
		\includegraphics[scale=0.30]{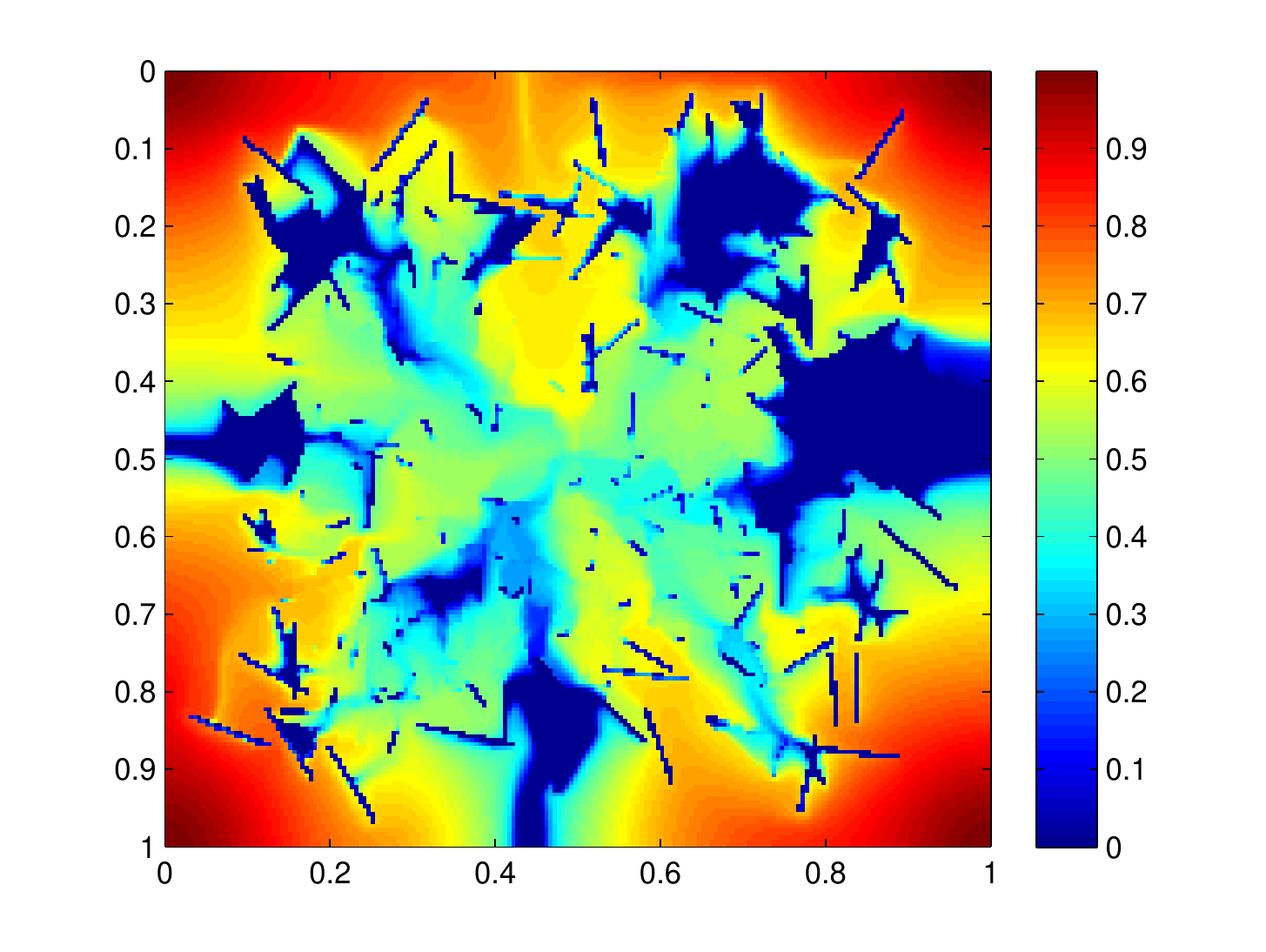}}	
	\subfigure[Ms-solution with 1 offline basis and 0 online basis, relative $L^2$ error is 39.5$\%$]{
		\includegraphics[scale=0.30]{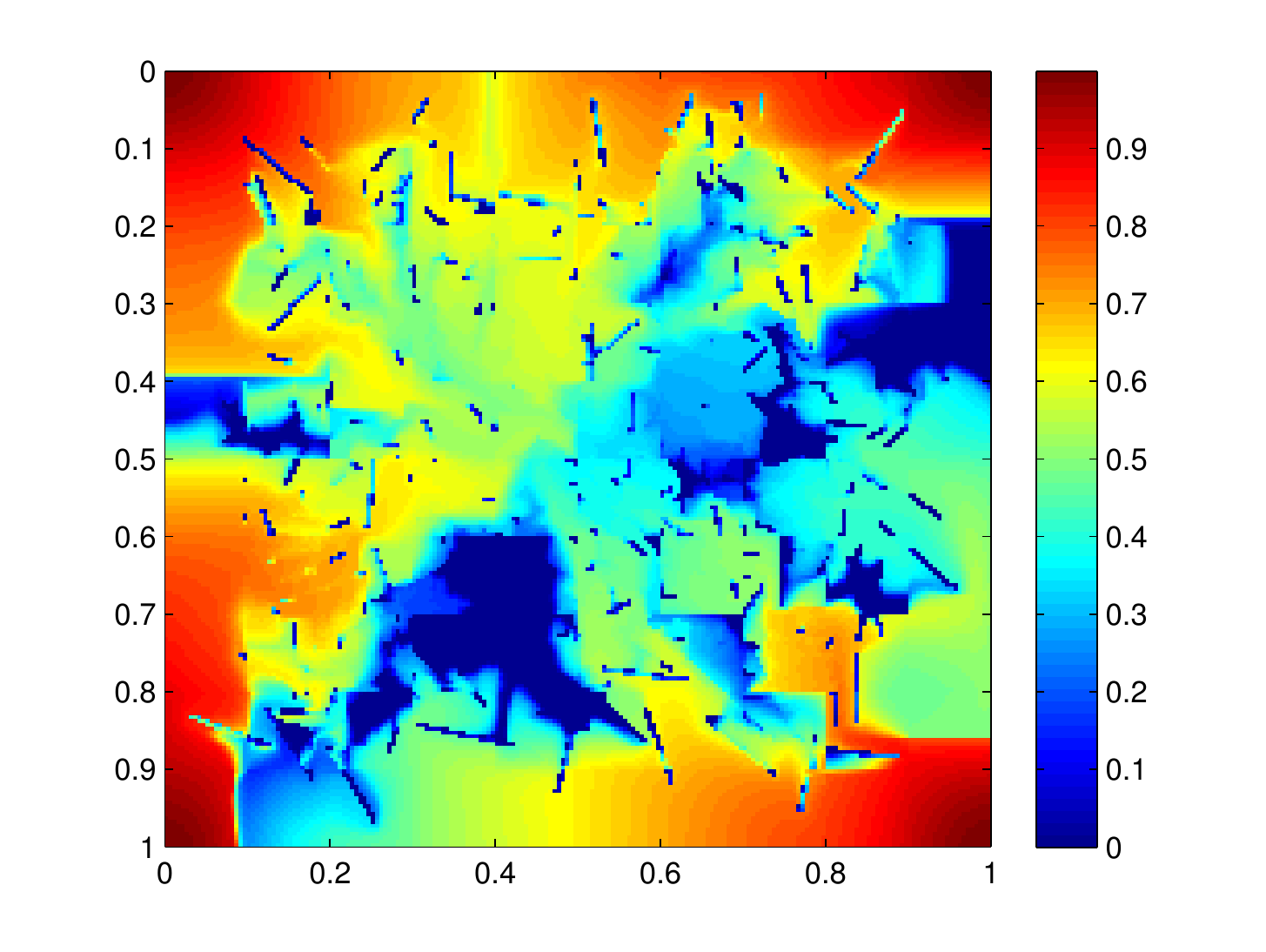}}	
	\subfigure[Ms-solution with 1 offline basis and 2 online basis, relative $L^2$ error is 8.2$\%$]{
		\includegraphics[scale=0.30]{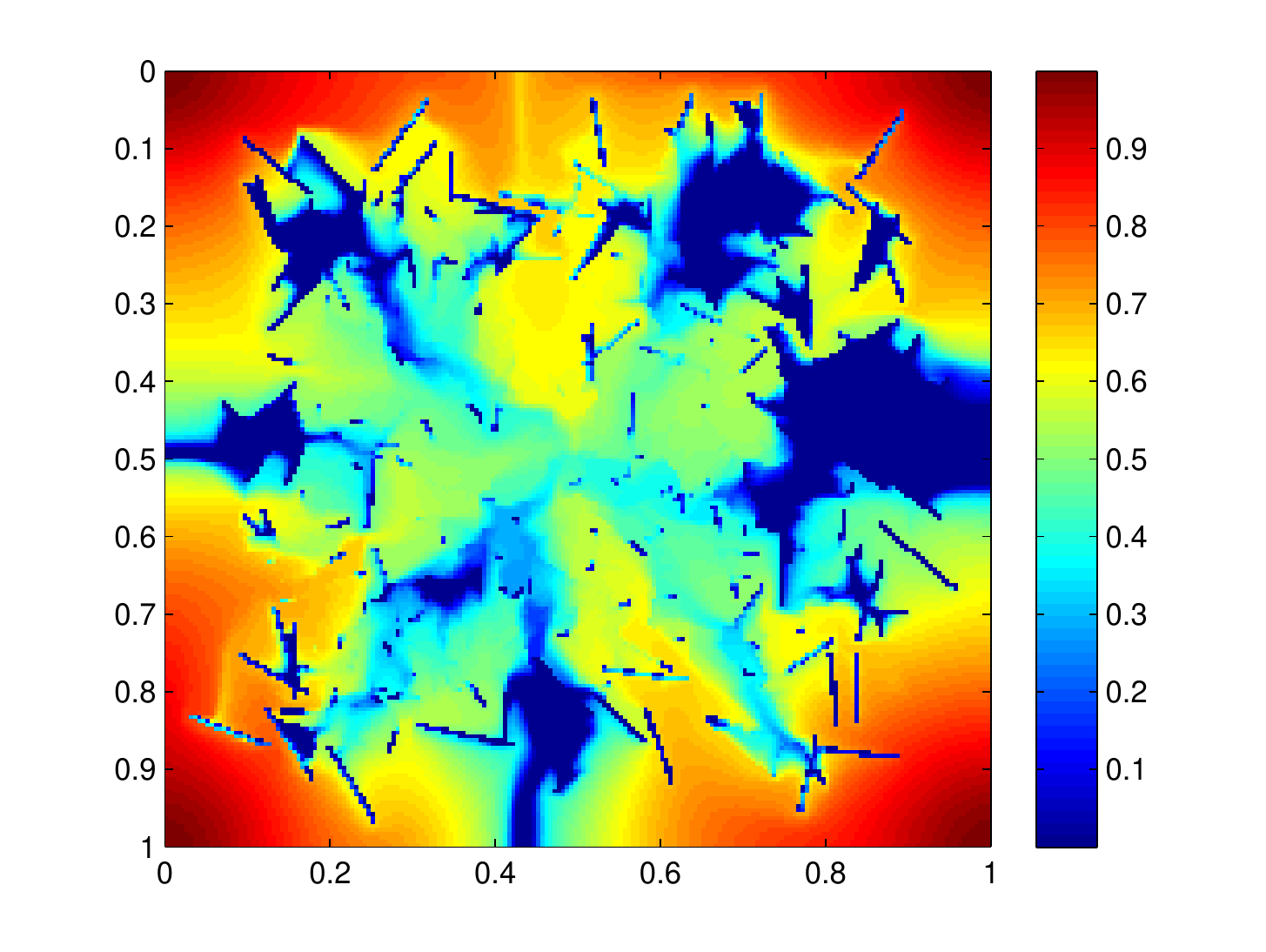}}

		
	\caption{Saturation comparison at $t=2000s$, model 1.}	
	\label{fig:Sat_modela}
\end{figure}

  Figure \ref{fig:Sat_modela} shows the saturation profile comparison (at time $t = 2000$) for model 1. We can see that if no online basis functions are
  used, the multiscale solution fails to capture many details  of the flow compared to the reference solution, which can be confirmed by the $39.5\%$ relative saturation error.
  If we use two online basis functions, the error drops to $8.2\%$, and the saturation profile
  is quite close to the reference profile.


\begin{figure}[H]
	\centering
	\subfigure{
		\includegraphics[scale=0.30]{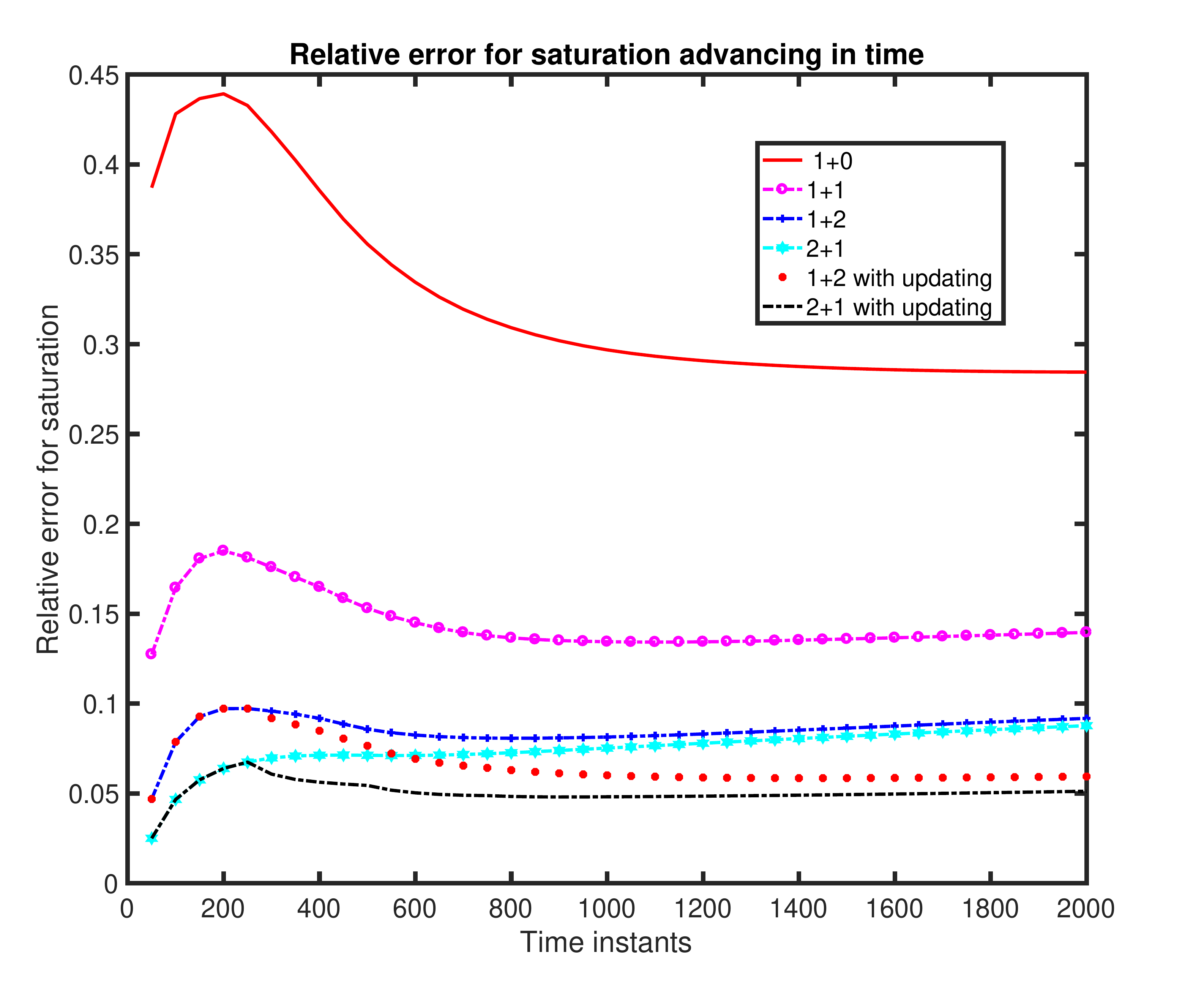}}	
		
	\caption{Saturation error for model 2.}	
	\label{fig:Sat_err50layers}
\end{figure}

\begin{figure}[H]
	\centering
	\subfigure[1 ]{
		\includegraphics[scale=0.2]{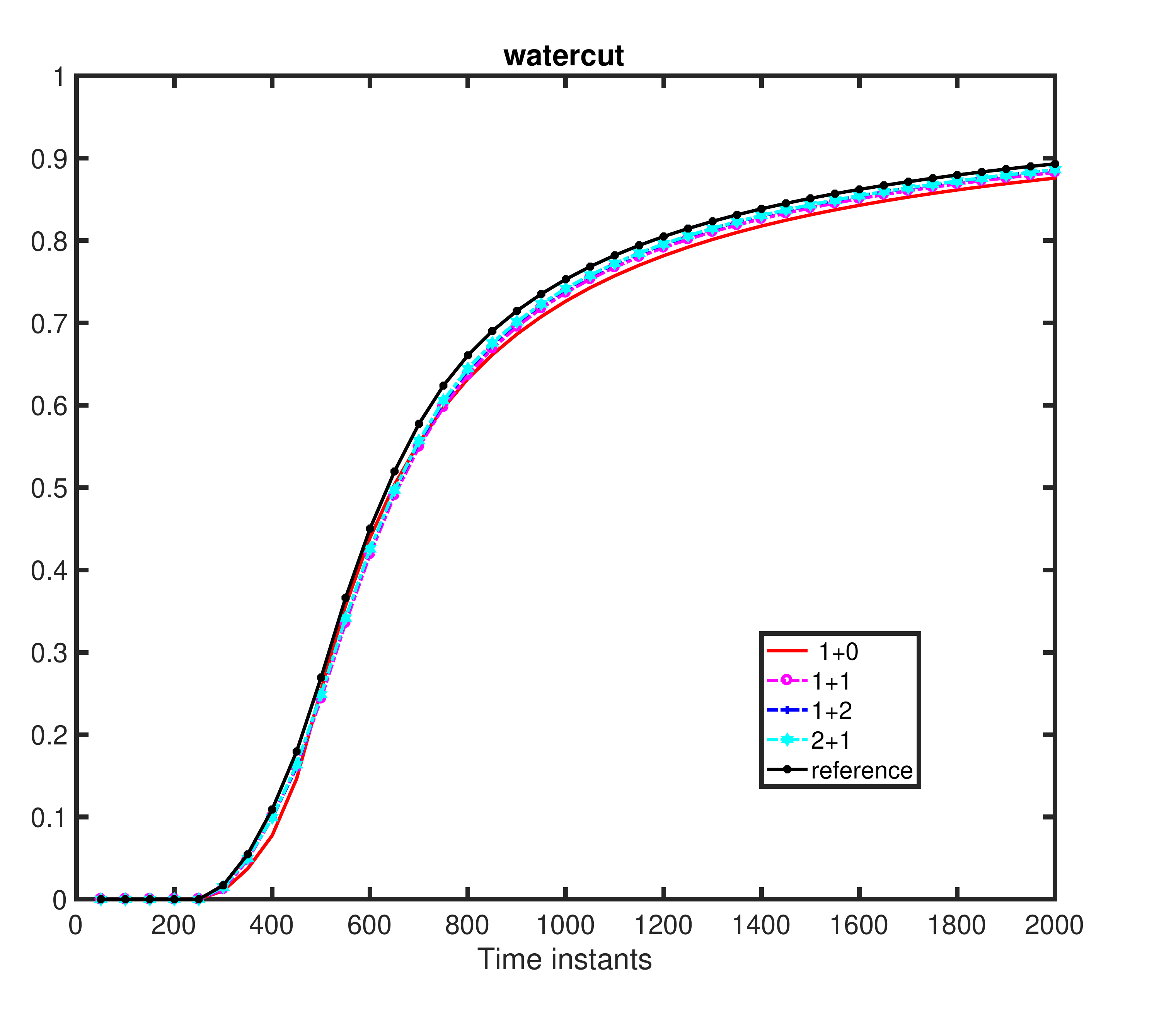}}	
	\subfigure[2]{
		\includegraphics[scale=0.20]{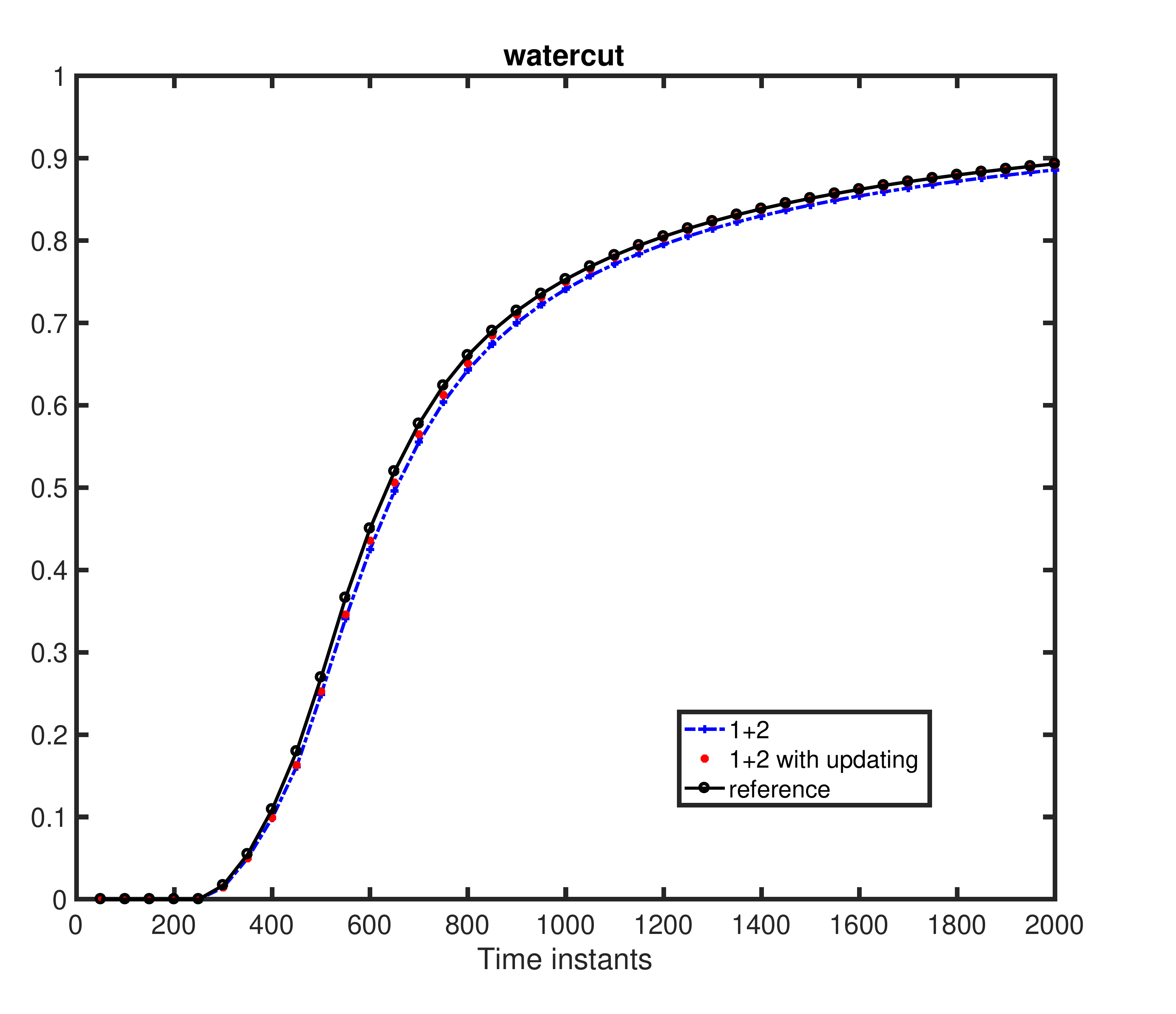}}	
	\subfigure[3]{
		\includegraphics[scale=0.20]{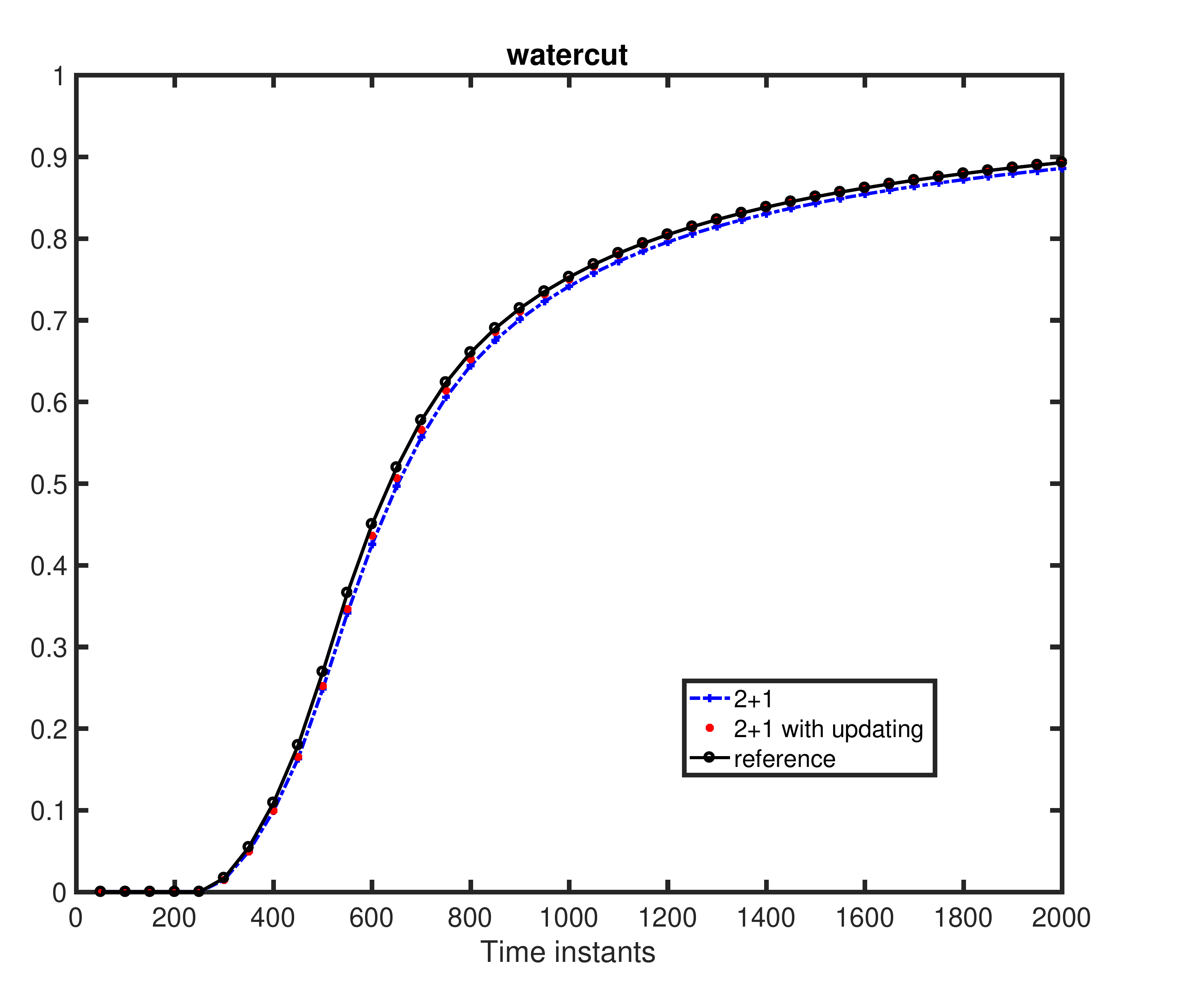}}

		
	\caption{Water-cut for model 2.}	
	\label{fig:watercut}
\end{figure}

\begin{figure}[H]
	\centering
	\subfigure[Fine-scale solution]{
		\includegraphics[scale=0.3]{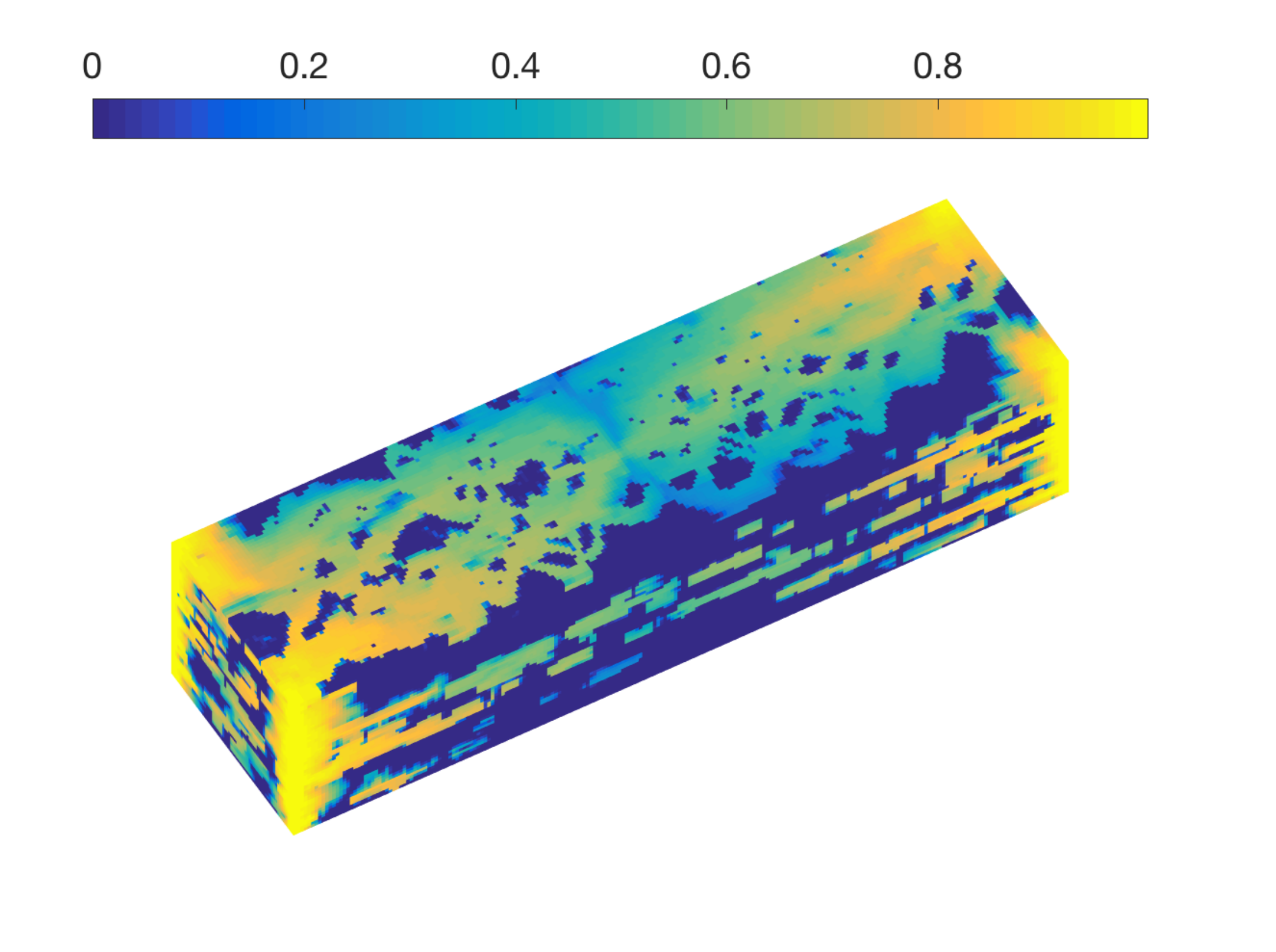}}	
	\subfigure[Ms-solution with 1 offline basis and 0 online basis, relative $L^2$ error is 28.4$\%$]{
		\includegraphics[scale=0.30]{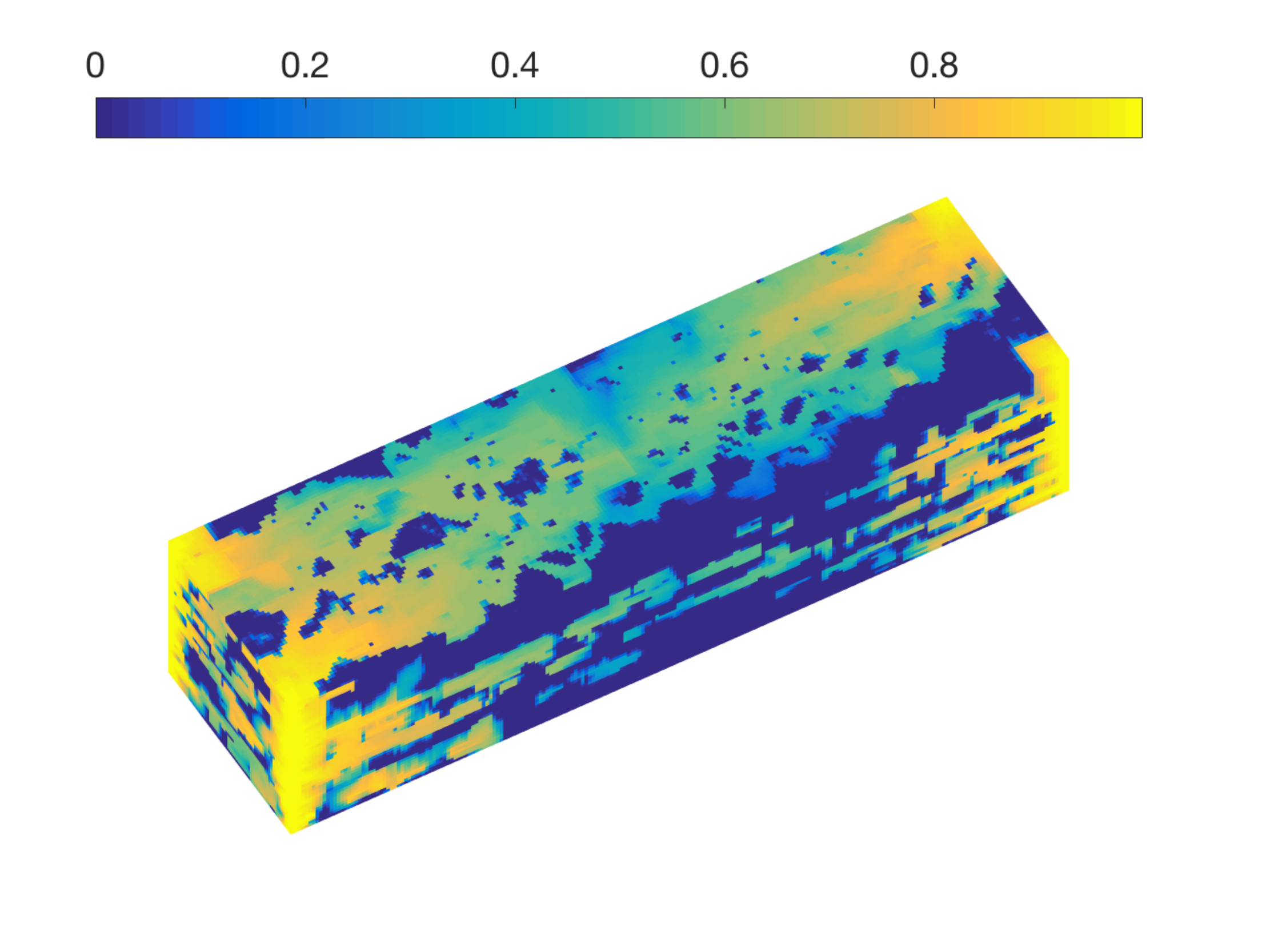}}	
	\subfigure[Ms-solution with 1 offline basis and 2 online basis, relative $L^2$ error is 9.2$\%$]{
		\includegraphics[scale=0.30]{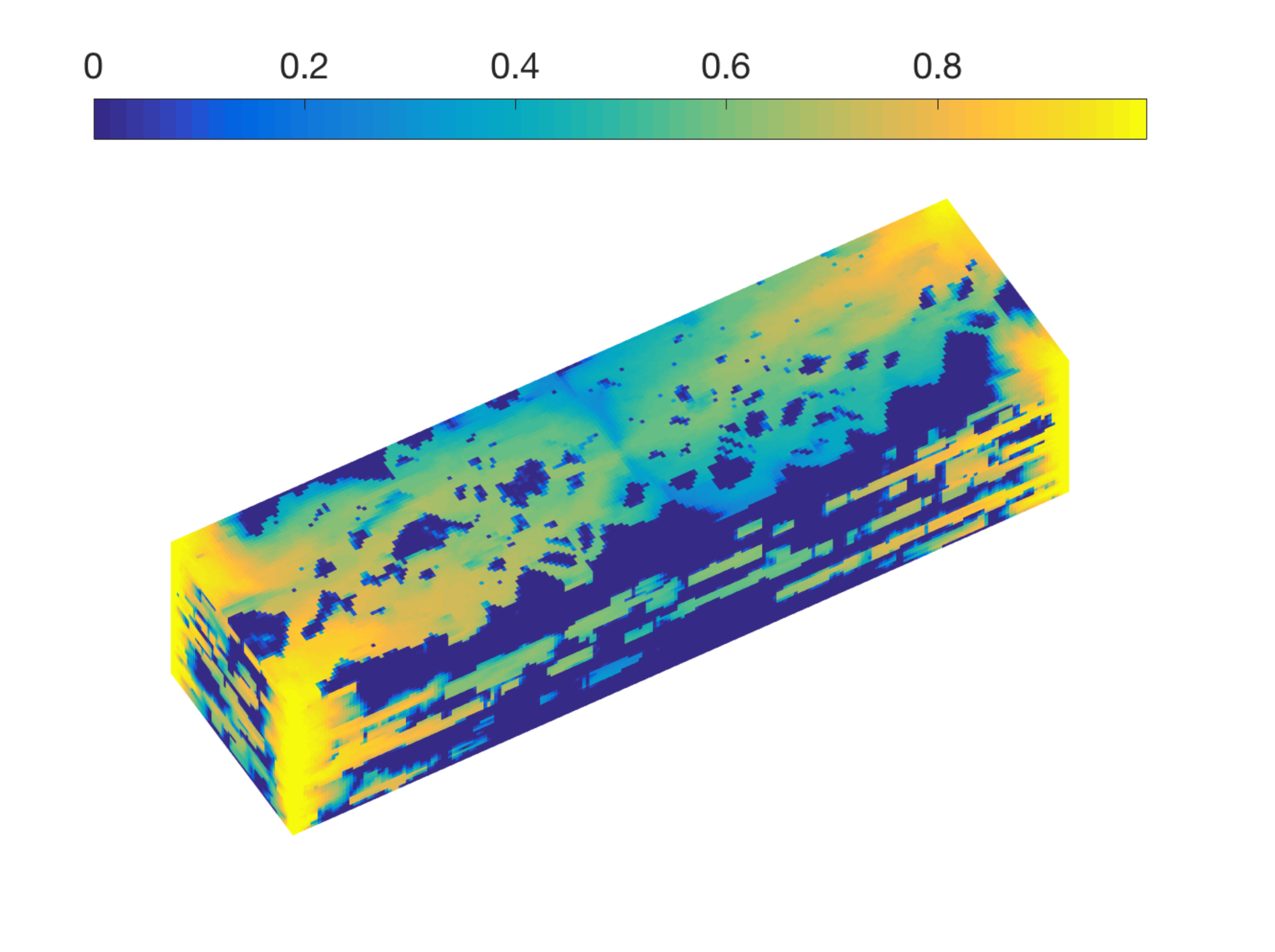}}
%
        \subfigure[Ms-solution with 1 offline basis and 2 online basis with basis updating every 5 time steps, relative $L^2$ error is 5.9$\%$]{
		\includegraphics[scale=0.30]{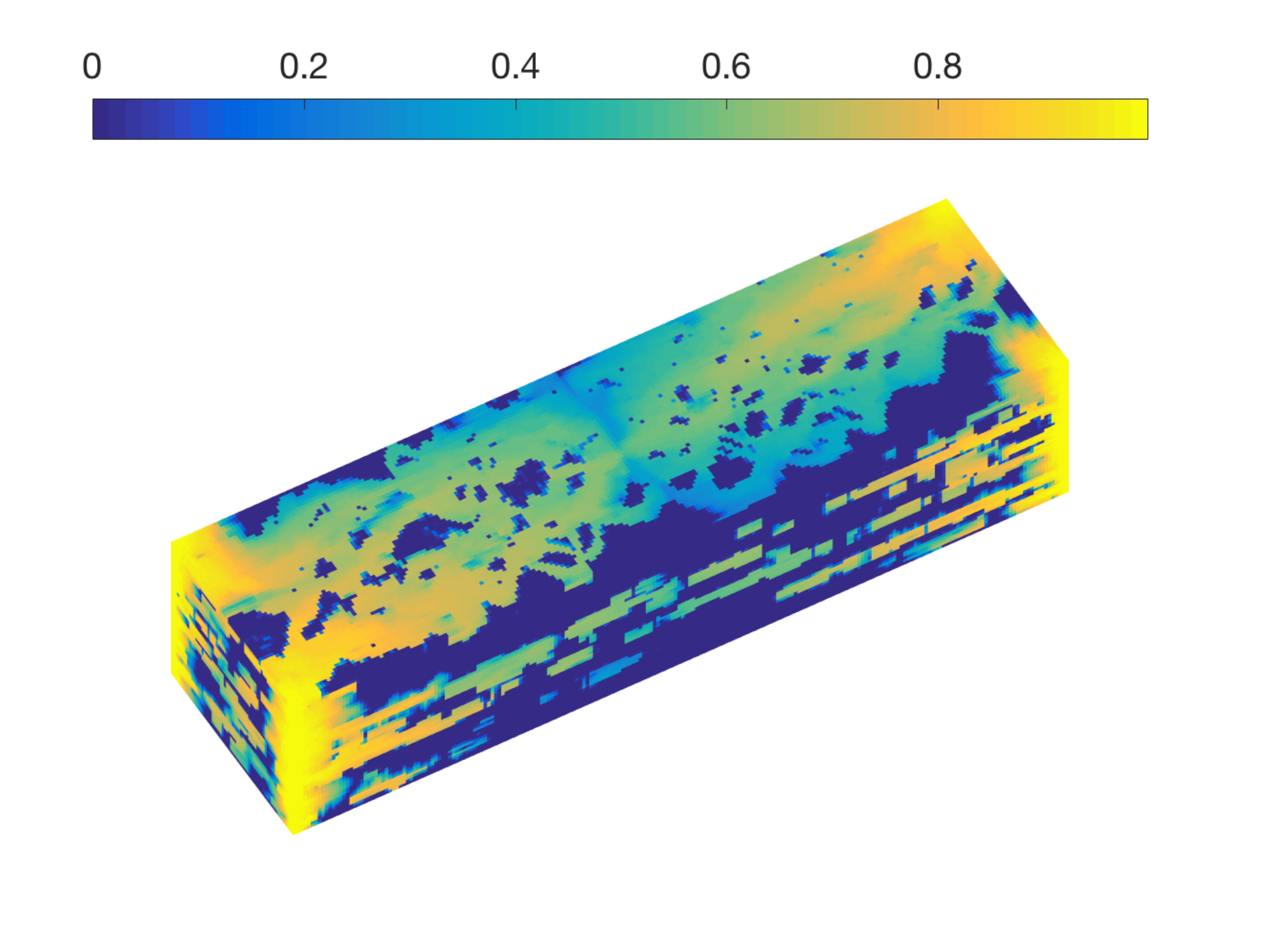}}
	\caption{Saturation for model 2}	
	\label{fig:saturation_3D}
\end{figure}

  We also obtain similar results for model 2. The relative error against
  time is shown in Figure \ref{fig:Sat_err50layers}. We also find that
  online basis is very efficient in reducing the saturation error.
  The water-cut comparison and saturation comparison against reference solution
  are show in Figure \ref{fig:watercut} and  \ref{fig:saturation_3D}, respectively. From Figure \ref{fig:watercut}, we can see that by using more online basis, the water cut gets closer to the reference one.  Updating the velocity basis at some middle time instants further improves the accuracy. In Figure \ref{fig:saturation_3D}, the saturation plots are from time $t=2000$.  The relative saturation error reduces from $28.4\%$ to $9.2\%$ by adding 2 online basis functions, and further reduces to $5.9\%$ by updating basis functions.

\section{Conclusions}
We have developed an oversampled online mixed generalized multiscale finite element method for flow problems in heterogeneous porous media. We start with an offline multiscale  solution and use the residual to construct online
basis  iteratively with oversampling techniques. The method produces mass conservative solutions for flow problems. By using
oversampled domain, we can compute online basis on a much smaller domain than a standard neighborhood of a coarse face, while gaining more accuracy.
Numerical results show that online basis with oversampling yields faster convergence speed compared with
no oversampling method. We also applied our method to flow and transport problems and demonstrate that by adding a few online basis functions, we can get good approximation solution.

\section*{Acknowledgments}

The research of Eric Chung is partially supported by RGC and CUHK.

	\bibliographystyle{plain}
	\bibliography{references}

\end{document}